\documentclass[11pt,a4paper]{article}
\usepackage[T2A]{fontenc}
\usepackage[cp1251]{inputenc}
\usepackage[russian]{babel}
\usepackage{array}
\usepackage{amsmath}
\usepackage{amsthm}
\usepackage{amssymb}
\usepackage{graphicx}
\usepackage{textcomp}

\providecommand{\keywords}{\textit{Keywords: }}
\newcommand{\sn}{\textrm sn}
\newcommand{\cn}{\textrm cn}
\newcommand{\dn}{\textrm dn}
\newcommand{\zn}{\textrm zn}

\begin{document}
\title{	On two-frequency quasi-periodic perturbations of systems close to two-dimensional Hamiltonian ones with a double limit cycle}
\author{O. S. Kostromina}
\date{{\small\it{Department of Differential Equations, Mathematical and Numerical Analysis,\\ Lobachevsky State University of Nizhny Novgorod, \\ 
	prospekt Gagarina, 23, Nizhny Novgorod,  603950, Russia \\
	os.kostromina@yandex.ru}}}
\maketitle
\maketitle MSC2010: 34C15
\renewcommand{\abstractname}{ }
\begin{abstract}
	The problem of the effect of two-frequency quasi-periodic perturbations on systems close to arbitrary nonlinear two-dimensional Hamiltonian ones is studied in the case when the corresponding perturbed autonomous systems have a double limit cycle. Its solution is important both for the theory of synchronization of nonlinear oscillations and for the theory of bifurcations of dynamical systems. In the case of commensurability of the natural frequency of the unperturbed system with frequencies of quasi-periodic perturbation, resonance occurs. Averaged systems are derived that make it possible to ascertain the structure of the resonance zone, that is, to describe the behavior of solutions in the neighborhood of individual resonance levels. The study of these systems allows determining possible bifurcations arising when the resonance level deviates from the level of the unperturbed system, which generates a double limit cycle in a perturbed autonomous system. The theoretical results obtained are applied in the study of a two-frequency quasi-periodic perturbed pendulum-type equation and are illustrated by numerical computations.
\end{abstract}


\section{Introduction}

In recent years there has been considerable interest in the study of multi-frequency quasi-periodic oscillations \cite{1}--\cite{13}. Such phenomena are widespread in nature and technology, in almost all areas of physics (radiophysics, electronics, laser physics, astrophysics, and many others), as well as in biology, chemistry, and medicine (see the above papers as well as monographs \cite{14,15,16}). Systems with a quasi-periodic dynamics are primarily non-autonomous system with two or more incommensurable frequencies. However, physicists have revealed the presence of quasi-periodic oscillations in autonomous systems (generators of quasi-periodic oscillations) \cite{17,18,19}.

Problems of the existence of quasi-periodic and almost periodic solutions, as well as complex dynamics in systems close to Hamiltonian ones, were studied in \cite{6}, \cite{20}--\cite{26}. The papers \cite{20,21} discussed the existence of quasi-periodic and almost periodic solutions for the Duffing equation. The questions of the existence of complex dynamics for the pendulum-type equation and the Duffing equation under quasi-periodic perturbations were considered in \cite{6}, \cite{22}--\cite{26}. More general results (for ordinary differential equations that contain a small parameter) on the existence of almost periodic solutions and integral manifolds are presented in \cite{27}. In contrast to these papers, papers \cite{3,7,8,11} are devoted to the problem of the effect of small nonconservative quasi-periodic perturbations on systems close to arbitrary nonlinear two-dimensional Hamiltonian ones with limit cycles. Based on research methods periodically perturbed systems  \cite{28,29}, the structure of non-degenerate resonance zones is revealed, conditions for the existence of quasi-periodic solutions in them are found, the problem of synchronization of quasi-periodic oscillations is solved, and the global behavior of solutions is analyzed. The theoretical results obtained in these studies are illustrated by the example of Duffing-type equations.

Let us consider the system
\begin{equation}\label{eq1}
\left\{\begin{aligned}
\dot{x}&=\frac{\partial H(x,y)}{\partial y}+\varepsilon g(x,y,\omega_1t, \omega_2t),\\
\dot{y}&=-\frac{\partial H(x,y)}{\partial x}+\varepsilon f(x,y,\omega_1t, \omega_2t),
\end{aligned}\right.
\end{equation}
where $\varepsilon$ is a small positive parameter, Hamiltonian function H and the functions $g$, $f$ are sufficiently smooth (analytic) 
and uniformly bounded in $x, y$ within some domain $D\subset\mathbb R^2$ (or $D\subset \mathbb R^1\times\mathbb S^1$) along with partial derivatives up to the second order. In addition, functions $g$, $f$ are continuous and quasi-periodic in $t$ uniformly with respect to $(x,y)\in D$ with incommensurable frequencies $\omega_1, \omega_2$. 

Let us assume that the unperturbed system ($\varepsilon=0$) is nonlinear and has a cell $D_0\in D$ filled with closed phase curves $H(x, y) = h$, $h\in [h_{min}, h_{max}]$, and not containing any equilibria and separatrices.
Let us suppose also that the perturbation is non-conservative, which is equivalent to the condition: $g'_x+f'_y \not\equiv 0$.

Passing from variables $x$, $y$ to the action $I$ -- angle $\theta$ variables in $D_0$:
\begin{equation}\label{eq2}
x=X(I,\theta), \ y=Y(I,\theta),
\end{equation}
where $X$, $Y$ are functions periodic in $\theta$ with period $2\pi$, we obtain the system of the form
\begin{equation}\label{eq3}
\left\{\begin{aligned}
\dot I&=\varepsilon F(I,\theta,\theta_1, \theta_2),\\
\dot\theta&=\omega(I)+\varepsilon G(I,\theta,\theta_1, \theta_2),\\
\dot\theta_1&=\omega_1,\\
\dot\theta_2&=\omega_2,
\end{aligned}\right.
\end{equation}
where $F(I,\theta,\theta_1, \theta_2)\equiv f(X,Y,\theta_1, \theta_2)X'_\theta-g(X,Y,\theta_1, \theta_2)Y'_\theta$, $G(I,\theta,\theta_1, \theta_2)\equiv -f(X,Y,\theta_1, \theta_2)X'_I+g(X,Y,\theta_1, \theta_2)Y'_I.$ Here $\omega(I)$ is the frequency of motion in closed phase curves (natural frequency of the unperturbed system). We assume that the function $\omega(I)$ is monotonic and does not vanish on the interval $(I_{min}, I_{max})\equiv (I(h_{min}), I(h_{max}))$. The functions $F, G$ are smooth enough in $I,\theta,\theta_1, \theta_2$ within domain $[I_{min}, I_{max}]\times \mathbb T^3$, where $\mathbb T^3$ is a three-dimensional torus. 

It is said that in system~\eqref{eq3} there is a \textit{resonance} if:
\begin{equation}\label{eq4}
n\omega(I)=m_1\omega_1+m_2\omega_2,
\end{equation}
where $n$, $m_1$ and $m_2$ are coprime integer numbers. We denote the value of $I$, calculated from (\ref{eq4}),
by $I_{nm_1m_2}$. Accordingly, we will refer to the unperturbed levels $I=I_{nm_1m_2}$ (closed phase curve $H(x,y)=h_{nm_1m_2}$ of the unperturbed system) as the \textit{resonance levels}.

Along with the system~\eqref{eq1}, let us consider the corresponding autonomous system
\begin{equation}\label{eq5}
\left\{\begin{aligned}
\dot{x}&=\frac{\partial H(x,y)}{\partial y}+\varepsilon g_0(x,y),\\
\dot{y}&=-\frac{\partial H(x,y)}{\partial x}+\varepsilon f_0(x,y),
\end{aligned}\right.
\end{equation}
where
\begin{equation*}
g_0(x,y)=\frac{1}{4\pi^2}\int\limits_{0}^{2\pi}\int\limits_{0}^{2\pi}g(x,y,\theta_1, \theta_2)d\theta_1 d\theta_2, \ f_0(x,y)=\frac{1}{4\pi^2}\int\limits_{0}^{2\pi}\int\limits_{0}^{2\pi}f(x,y,\theta_1, \theta_2)d\theta_1 d\theta_2.
\end{equation*}

Suppose that system~\eqref{eq5} has a double limit cycle. This means that the Poincar\'e--Pontryagin generating function has a double root $h=h_{*}$ (for more details, see the next section).

The effect of quasi-periodic perturbation for a system of the form (\ref{eq5}) with a structurally stable limit cycle (in this case, the Poincar\'e--Pontryagin generating function has a simple root $h=h_0$) was studied in \cite{1,3,7,8,11}. The action of time-periodic perturbation for a system with a double limit cycle was investigated in \cite{30}. However, this study excluded the case of parametric perturbations. A Duffing-type equation was considered as an illustrative example. 

The problem of the effect of quasi-periodic perturbations on a perturbed autonomous system of the form (\ref{eq5}), which has a double limit cycle, has not been previously studied and is considered for the first time in this paper. We study the rearrangements of phase portraits of the averaged system near the bifurcation case, when the resonance level $I=I_{nm_1m_2}$ coincides with the level $I=I_*$ in the neighborhood of which system (\ref{eq5}) has a double limit cycle. The solution to this problem is essential for the theory of synchronization of oscillations, as well as the theory of bifurcations of dynamical systems. This study largely follows the paper \cite{30}, as well as the works \cite{3,7,11}.

As an illustrative example, we consider a pendulum-type equation of the form
\begin{equation}\label{eqExample}
\ddot{x}+\sin{x}=\varepsilon [(-1+p_1\cos{3x}+p_2x\alpha (t))\dot{x}+p_3\alpha (t)],
\end{equation}
where $p_1$, $p_2>0$, $p_3>0$ are parameters, $\varepsilon$ is a small positive parameter, $\alpha (t)=\cos{\omega_1t}\sin{\omega_2t}$. We suppose that $\omega_1$, $\omega_2$ are incommensurable. The perturbation contains a nonlinear parametric term $p_2x\dot{x}\alpha (t)$. The case when $p_2=0$ and the function $\alpha (t)$ is periodic was studied in~\cite{31}.

\section{Poincar\'e--Pontryagin generating function and double cycles in the perturbed autonomous system}

In cell $D_0$, pass to the action $I$ -- angle $\theta$ variables in system~\eqref{eq5}, writing this transformation in the form~\eqref{eq2}, and averaging the resulting system over $\theta$, we obtain the system
\begin{equation}\label{eq6}
\dot u=\varepsilon B_0(u),
\end{equation}
where
\begin{equation}\label{eq7}
B_0(u)=\frac 1{2\pi}\int\limits_0^{2\pi}[f_0(X,Y)X'_\theta-g_0(X,Y)Y'_\theta]\,d\theta.
\end{equation}
Here, $u=I+O(\varepsilon)$. The function $B_0(u)$ is called the Poincar\'e--Pontryagin generating function.

It is known \cite{28,29}, if $u=u_0$ is a simple equilibrium state ($B_0(u_0)=0$, $B'_0(u_0)\neq 0$) of the averaged system~\eqref{eq6}, then in system~\eqref{eq5} for sufficiently small $\varepsilon$ this equilibrium state corresponds to a rough limit cycle. Moreover, the limit cycle will be stable if the equilibrium state is stable, that is, $B'(u_0)<0$, and unstable otherwise. Let us denote
\begin{equation*}
B_1(u_0)\equiv B'_0(u_0)=\frac 1{2\pi}\int\limits^{2\pi }_{0}\left({g_0}'_x+{f_0}'_y\right)_{\left|{\begin{subarray}{c}x=X(u_0,\theta)\\
		y=Y(u_0,\theta)\end{subarray} }\right.}d\theta.
\end{equation*}

The double root $u=u_*$ of the equation $B_0(u)=0$ determines the unperturbed level $I=I_*$ (closed phase curve $h=h_*$ of the unperturbed system), from which a double limit cycle is generated under the action of the perturbation. In this case, $B_1(u_*)=0$, $B'_1(u_*)\neq 0$.

When non-autonomous disturbances affect an autonomous system, a limit cycle may coincide with the resonance level.

\section{Averaging system in the neighborhood of the individual resonance level}

It is known~\cite{3,7}, in the neighborhood
$U_\mu=\{(I, \theta ): I_{nm_1m_2}-C\mu<I<I_{nm_1m_2}+C\mu, \ 0\leqslant \theta <2\pi, \  C=\textrm{const}>0 \}, \ \mu=\sqrt\varepsilon $ of the individual resonance level $I = I_{nm_1m_2}$ (it will be called the \textit{resonance zone}) the system~\eqref{eq2}  reduces to an averaged system of the form
\begin{equation}\label{eq8}
\left\{\begin{aligned}
\dot{u} &=\mu A(v;I_{nm_1m_2})+\mu^2P_0(v;I_{nm_1m_2})u+O(\mu^3),\\
\dot{v} &=\mu b_1u+\mu^2(b_2u^2+Q_0(v;I_{nm_1m_2}))+O(\mu^3),
\end{aligned}\right.
\end{equation}
where  
\begin{equation*}
A(v;I_{nm_1m_2})=\frac 1{4\pi^2 n^2}\int\limits_0^{2\pi n}\int\limits_0^{2\pi n}\left.F(I_{nm_1m_2}, v+\frac{m_1\theta_1+m_2\theta_2}{n}, \theta_1,\theta_2) \right.\,d\theta_1 d\theta_2 ,
\end{equation*}
\begin{equation*}
P_0(v;I_{nm_1m_2})=\frac 1{4\pi^2 n^2}\int\limits_0^{2\pi n}\int\limits_0^{2\pi n}\left.F'_I(I_{nm_1m_2}, v+\frac{m_1\theta_1+m_2\theta_2}{n}, \theta_1,\theta_2) \right.\,d\theta_1 d\theta_2 ,
\end{equation*}
\begin{equation*}
Q_0(v;I_{nm_1m_2})=\frac 1{4\pi^2 n^2}\int\limits_0^{2\pi n}\int\limits_0^{2\pi n}\left.G(I_{nm_1m_2}, v+\frac{m_1\theta_1+m_2\theta_2}{n}, \theta_1,\theta_2) \right.\,d\theta_1 d\theta_2 ,
\end{equation*}
\begin{equation*}
b_1= \omega'(I_{nm_1m_2}), \ b_2= \omega''(I_{nm_1m_2})/2.
\end{equation*}

The functions $A(v;I_{nm_1m_2})$, $P_0(v;I_{nm_1m_2})$ and $Q_0(v;I_{nm_1m_2})$ are periodic in $v$ with the least period equal to $2\pi/n$~\cite{3,7}. Therefore, the phase space of system~\eqref{eq8} is the cylinder $\{v(mod(2\pi/n)),u\}$. A simple stable (unstable) equilibrium state of the system~\eqref{eq8} corresponds to a stable (unstable) quasi-periodic resonance solution with periods $2\pi n/\omega_1$, $2\pi n/\omega_2$, and two-dimensional invariant tori with quasi-periodic motion in the original four-dimensional system.

The functions $A(v;I_{nm_1m_2})$ and $P_0(v;I_{nm_1m_2})$ can be represented as
\begin{equation*}
\begin{array}{cc}
A(v;I_{nm_1m_2})=\widetilde{A}(v;I_{nm_1m_2})+B_0(I_{nm_1m_2}), P_0(v;I_{nm_1m_2})=\widetilde{P}_0(v;I_{nm_1m_2})+B_1(I_{nm_1m_2}),\\
B_0(I_{nm_1m_2})=\left<A(v;I_{nm_1m_2})\right>_v=\displaystyle\frac{n}{2\pi}\int\limits_0^{2\pi/ n}A(v;I_{nm_1m_2})\,dv,\\ B_1(I_{nm_1m_2})=\left<P_0(v;I_{nm_1m_2})\right>_v=\displaystyle\frac{n}{2\pi}\int\limits_0^{2\pi/ n}P_0(v;I_{nm_1m_2})\,dv.
\end{array}
\end{equation*}

Since we consider the case, when $I_{nm_1m_2}=I_*$, then $B_0(I_{nm_1m_2})=B_1(I_{nm_1m_2}) = 0$, and therefore, it is necessary to take into account the terms $O(\mu^3)$.  Neglecting the terms $O(\mu^4)$, we obtain the system 
\begin{equation}\label{eq9}
\left\{\begin{aligned}
\frac{du}{d\tau} &= \widetilde{A}(v;I_{nm_1m_2})+B_0(I_{nm_1m_2})+\mu\left[\widetilde{\sigma}(v;I_{nm_1m_2})+ B_1(I_{nm_1m_2})\right]+\mu^2\left[\left(\widetilde{P}_1(v;I_{nm_1m_2})\right.\right.+\\&+\left.\left.B_2(I_{nm_1m_2})+\frac{b_2}{b_1}\frac{dQ_0(v;I_{nm_1m_2})}{dv}\right)u^2-\widetilde{P}_0(v;I_{nm_1m_2})\frac{Q_0(v;I_{nm_1m_2})}{b_1}\right],\\
\frac{dv}{d\tau} &= b_1u+\mu b_2u^2+\mu^2\left[b_3u^3+\left(Q_1(v;I_{nm_1m_2})-2\frac{b_2}{b_1}Q_0(v;I_{nm_1m_2})\right)u \right],
\end{aligned}\right.
\end{equation}
where 
\begin{equation*}
P_1(v;I_{nm_1m_2})=\frac 1{8\pi^2 n^2}\int\limits_0^{2\pi n}\int\limits_0^{2\pi n}\left.F''_{II}(I_{nm_1m_2}, v+\frac{m_1\theta_1+m_2\theta_2}{n}, \theta_1,\theta_2) \right.\,d\theta_1 d\theta_2 ,
\end{equation*}
\begin{equation*}
Q_1(v;I_{nm_1m_2})=\frac 1{4\pi^2 n^2}\int\limits_0^{2\pi n}\int\limits_0^{2\pi n}\left.G'_I(I_{nm_1m_2}, v+\frac{m_1\theta_1+m_2\theta_2}{n}, \theta_1,\theta_2) \right.\,d\theta_1 d\theta_2 ,
\end{equation*}
\begin{equation*}
b_3= \omega'''(I_{nm_1m_2})/6.
\end{equation*}
The functions $P_1(v;I_{nm_1m_2})$ can be represented as $P_1(v;I_{nm_1m_2})=\widetilde{P}_1(v;I_{nm_1m_2})+B_2(I_{nm_1m_2})$, $B_2(I_{nm_1m_2})=B''_0(I_{nm_1m_2})/2$.

Making the change $u\to u - \mu Q_0(v, I_{nm_1m_2})/b_1$ in system~\eqref{eq9} and passing to the slow time $\tau=\mu t$, we obtain the system
\begin{equation}\label{eq10}
\left\{\begin{aligned}
\frac{du}{d\tau} &= \widetilde{A}(.)+B_0(I_{nm_1m_2})+\mu\left[\widetilde{\sigma}(.)+ B_1(I_{nm_1m_2})\right]+\mu^2\left[\left(\widetilde{P}_1(.)+B_2(I_{nm_1m_2})\right.\right.+\\&+\left.\left.\frac{b_2}{b_1}\frac{dQ_0(.)}{dv}\right)u^2-\widetilde{P}_0(.)\frac{Q_0(.)}{b_1}\right],\\
\frac{dv}{d\tau} &= b_1u+\mu b_2u^2+\mu^2\left[b_3u^3+\left(Q_1(.)-2\frac{b_2}{b_1}Q_0(.)\right)u \right],
\end{aligned}\right.
\end{equation}
where $(.)\equiv (v;I_{nm_1m_2})$, $\widetilde{\sigma}(v;I_{nm_1m_2})=\widetilde{P}_0(v;I_{nm_1m_2})+dQ_0(v;I_{nm_1m_2})/dv$.
This system determines the topology of individual resonance zones up to terms of order $\mu^3$. According to~\cite{3}, a simple stable (unstable) equilibrium state of the averaged system~\eqref{eq10} corresponds to a stable (unstable) quasi-periodic resonance solution of periods $\displaystyle\frac{2\pi n}{\omega_1}$, $\displaystyle\frac{2\pi n}{\omega_2}$ in system~\eqref{eq1}, or a two-dimensional stable (unstable) invariant torus in system~\eqref{eq3}. A rough limit cycle of the averaged system~\eqref{eq10} with frequency $\omega_0$ corresponds to a quasi-periodic resonance solution of periods $\displaystyle\frac{2\pi }{\omega_0}$, $\displaystyle\frac{2\pi n}{\omega_1}$, $\displaystyle\frac{2\pi n}{\omega_2}$ in system~\eqref{eq1}, or a three-dimensional invariant torus in system~\eqref{eq3} of the same stability type as the cycle.

An important role in the study of the structure of the neighborhood of resonance levels $I=I_{nm_1m_2}$ is played by the function $\sigma(v;I_{nm_1m_2})=\widetilde{\sigma}(v;I_{nm_1m_2})+B_1(I_{nm_1m_2})$, $B_1(I_{nm_1m_2})=\left<\sigma(v;I_{nm_1m_2})\right>_v$. First of all, it is necessary to find out whether it depends on $v$ or is constant~\cite{28,29}.

\section{Study of the averaged system}

Put $B_0(I_{nm_1m_2})=B_1(I_{nm_1m_2}) = 0$. Then system~\eqref{eq10} without taking into account the conservative term $2\mu\displaystyle\frac{b_2}{b_1}\widetilde{A}(v;I_{nm_1m_2})\dot{v}$ is transformed to the form
\begin{equation}\label{eq12}
\ddot{v}-b_1\widetilde{A}(v;I_{nm_1m_2})=\mu\widetilde{\sigma}(v;I_{nm_1m_2})\dot{v}+\mu^2\left[M(v;I_{nm_1m_2})\dot{v}^2+N(v;I_{nm_1m_2})\right],
\end{equation}
where
\begin{equation*}
M(v;I_{nm_1m_2})=\frac{1}{b_1}\left[\widetilde{P}_1(.)+{Q_1}'_v(.)+\frac{b_2}{b_1}\left(\widetilde{P}_0(.)+\widetilde{\sigma}(.)\right)+3\frac{b_3}{b_1}\widetilde{A}(.)\right]+\frac{1}{b_1}B_2(I_{nm_1m_2}),
\end{equation*}
\begin{equation*}
N(v;I_{nm_1m_2})=-\widetilde{P}_0(.)Q_0(.)+\left(Q_1(.)-2\frac{b_2}{b_1}Q_0(.)\right)\widetilde{A}(.).
\end{equation*}
Неre $(.)\equiv (v;I_{nm_1m_2})$.

Let us rewrite equation~\eqref{eq12} as 
\begin{equation}\label{eq13}
\left\{\begin{aligned}
\dot{u} &= \widetilde{A}(v;I_{nm_1m_2})+\mu\widetilde{\sigma}(v;I_{nm_1m_2})u+\frac{\mu^2}{b_1}\left[M(v;I_{nm_1m_2})b_1^2u^2+N(v;I_{nm_1m_2})\right],\\
\dot{v} &= b_1u
\end{aligned}\right. 
\end{equation}
Let us consider the case, when $\sigma(v;I_{nm_1m_2})=const$, therefore $\widetilde{\sigma} = 0$. In this case the system~\eqref{eq13} is invariant under the change $u\to -u, \tau\to -\tau$, therefore, the phase portrait of equation~\eqref{eq12} is symmetric about the $u_2$ axis. Phase space of system~\eqref{eq13} is the cylinder $\{v(mod(2\pi/n)),u\}$.
System~\eqref{eq13} is close to Hamiltonian with the Hamiltonian function $\widetilde{H}(v,u)=b_1u^2/2-V(v)$, $V(v)=\int\widetilde{A}(v;I_{nm_1m_2})dv$. Denote by $(v_0,0)$ the equilibrium state of the unperturbed system. System~\eqref{eq13} has simple equilibrium states of two types: the center, if $b_1\widetilde{A}'_v(v_0)<0$, and the saddle, if $b_1\widetilde{A}'_v(v_0)>0$. To establish the relative position of the separatrices of the saddle $(v_0,0)$ of the unperturbed system under the action of the perturbation, we use the Melnikov formula~\cite{32}: $\Delta_\mu=\mu^2\Delta_1+O(\mu^3)$, which determines (up to terms of order $\mu^3$) the distance between the perturbed separatrices. Making the substitution $v=w+v_0$, moving the equilibrium state $(v_0,0)$ of the saddle type of the unperturbed system to the origin, we obtain
\begin{equation}\label{eq14}
\Delta_1=\frac{1}{b_1}\int\limits^{\infty }_{-\infty}\left[M(w+v_0)b_1^2u^2+N(w+v_0)\right]\frac{dw}{d\tau}d\tau,
\end{equation}
where $w(\tau)$, $u(\tau)$ is the solution of the unperturbed system on the separatrix. From the integral of the unperturbed system we find
\begin{equation}\label{eq15}
u=\pm\sqrt{\frac{2}{b_1}\left(V(w,v_0)-V(0,v_0)\right)}, \ V(w,v_0)=\int\widetilde{A}(w+v_0)dw.
\end{equation}
Substituting~\eqref{eq15} into~\eqref{eq14}, we obtain
\begin{equation}\label{eq16}
\Delta_1=\frac{1}{b_1}\int\limits^{2\pi/n }_{0}\left[2b_1M(w+v_0)\left(V(w,v_0)-V(0,v_0)\right)+N(w+v_0)\right]dw.
\end{equation}

\begin{figure}[h!]
	\begin{center}
		\begin{tabular}{ccc}
			\includegraphics[width=140pt,height=140pt]{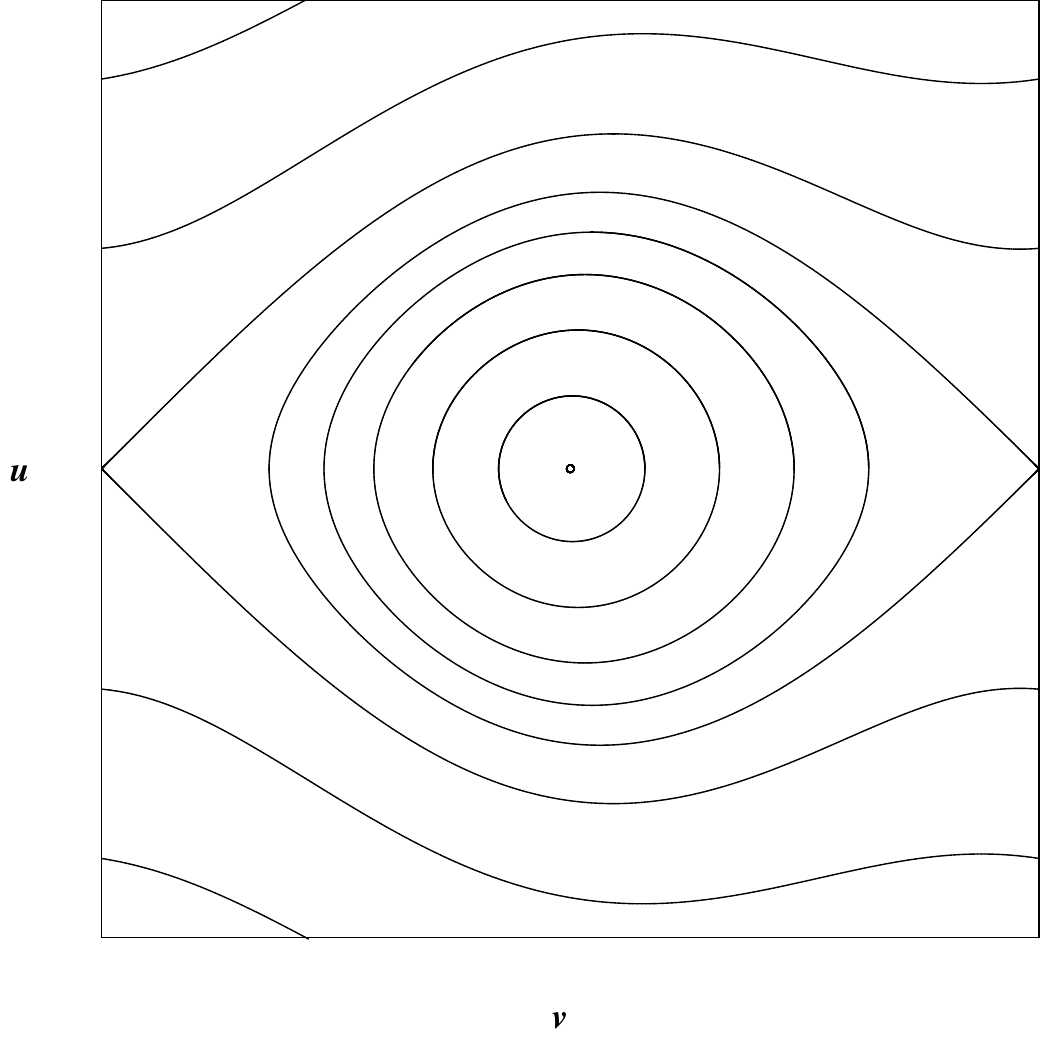}&
			\includegraphics[width=140pt,height=140pt]{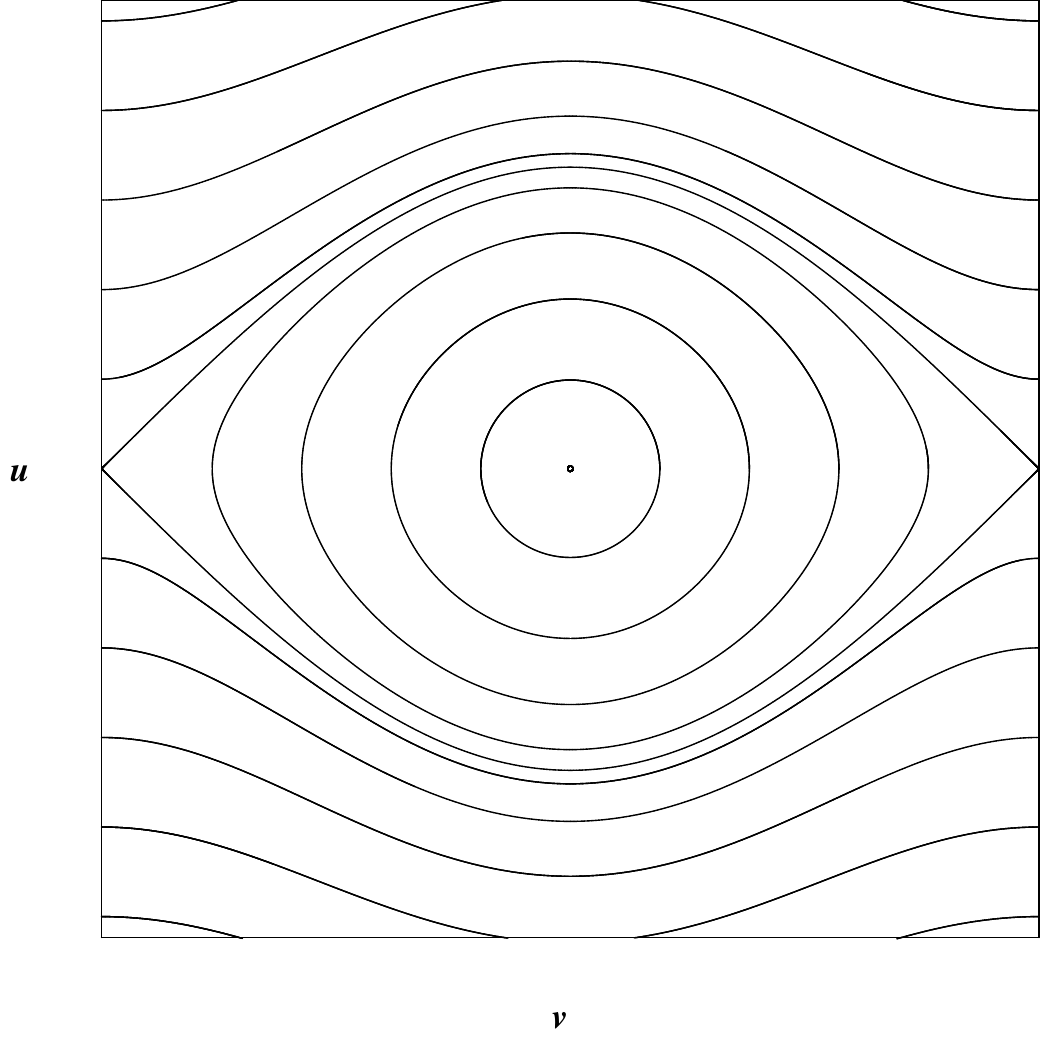}&
			\includegraphics[width=140pt,height=140pt]{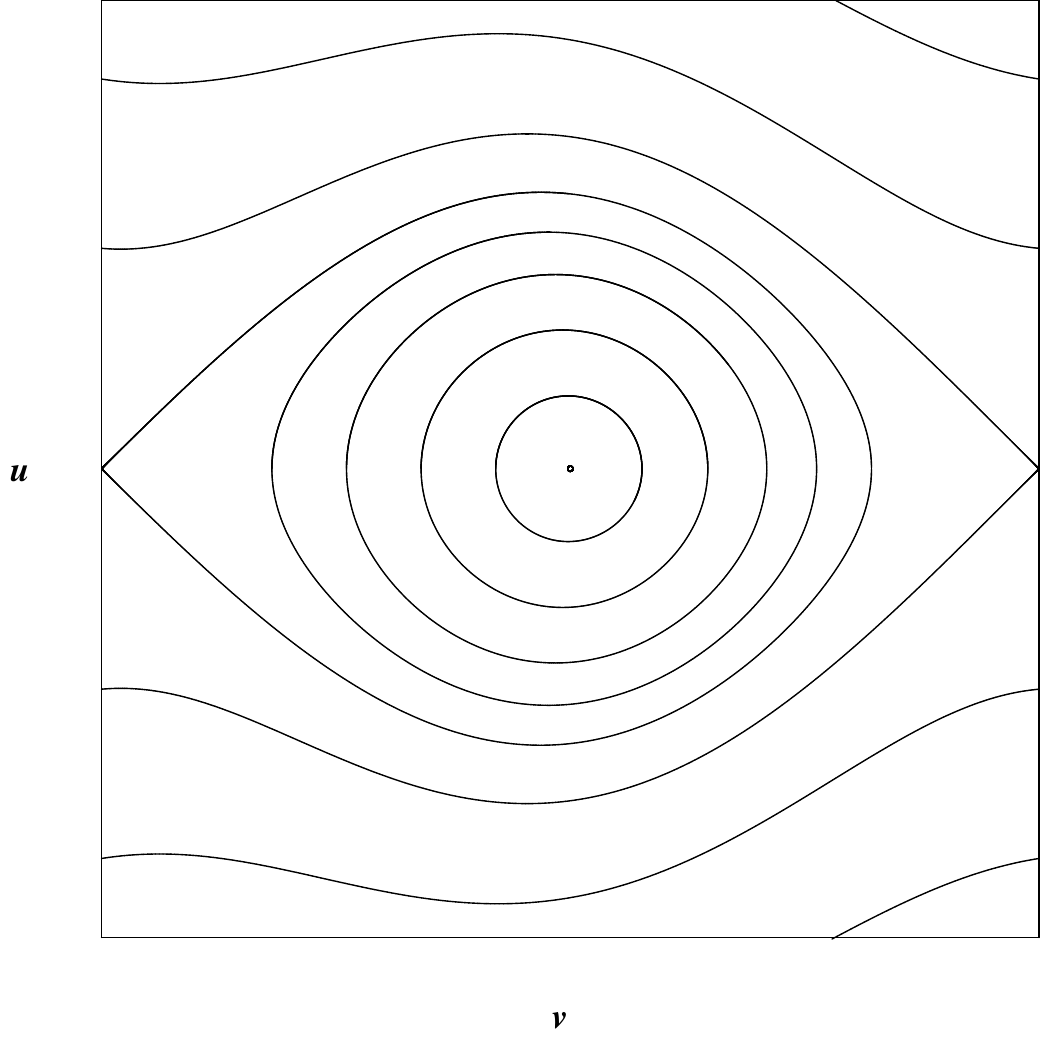}\\
			\footnotesize{(a) $B_2(I_{nm_1m_2})\neq 0$}& \footnotesize{(b) $B_2(I_{nm_1m_2})=0$}& \footnotesize{(c) $B_2(I_{nm_1m_2})\neq 0$}\\
		\end{tabular}
	\end{center}
	\renewcommand{\figurename}{Fig}
	\caption{Possible phase portraits of system~\eqref{eq13} for $\widetilde{\sigma} = 0$.}\label{fig1}
\end{figure}

There are two possible cases: $\Delta_1=0$ and $\Delta_1\neq 0$. In the first case, we have two separatrix loops of the saddle on the phase cylinder (on the upper and lower half-cylinders). This is a bifurcation case. In the second case, the separatrices are split, and the size of the gap depends on the relationship between the amplitude of the function  $\widetilde{M}(v;I_{nm_1m_2})=\widetilde{P}_1(v;I_{nm_1m_2})+{Q_1}'_v(v;I_{nm_1m_2})+\frac{b_2}{b_1}\widetilde{P}_0(v;I_{nm_1m_2})+3\frac{b_3}{b_1}\widetilde{A}(v;I_{nm_1m_2})$ and the value $B_2(I_{nm_1m_2})$. Possible phase portraits of system~\eqref{eq13} on the phase cylinder $\{v(mod(2\pi/n)),u\}$ are shown in Fig~\ref{fig1}.
The case of impassable resonance (synchronization of oscillations) is shown in Fig~\ref{fig1}(b): for any initial conditions of the phase space of the averaged system, the phase curves of the averaged system remain in a limited region of the phase space. The case of partly passable resonance is shown in Figs~\ref{fig1}(a) and (c): in the phase space of the averaged system, there are sets of initial conditions (inside the separatrix loop) that correspond to bounded phase curves, as well as sets of initial conditions that correspond to phase curves (outside the separatrix loop) that leave any bounded region of the phase space when $\tau\to\infty$.

\begin{figure}[h!]
	\begin{center}
		\begin{tabular}{cc}
			\includegraphics[width=140pt,height=140pt]{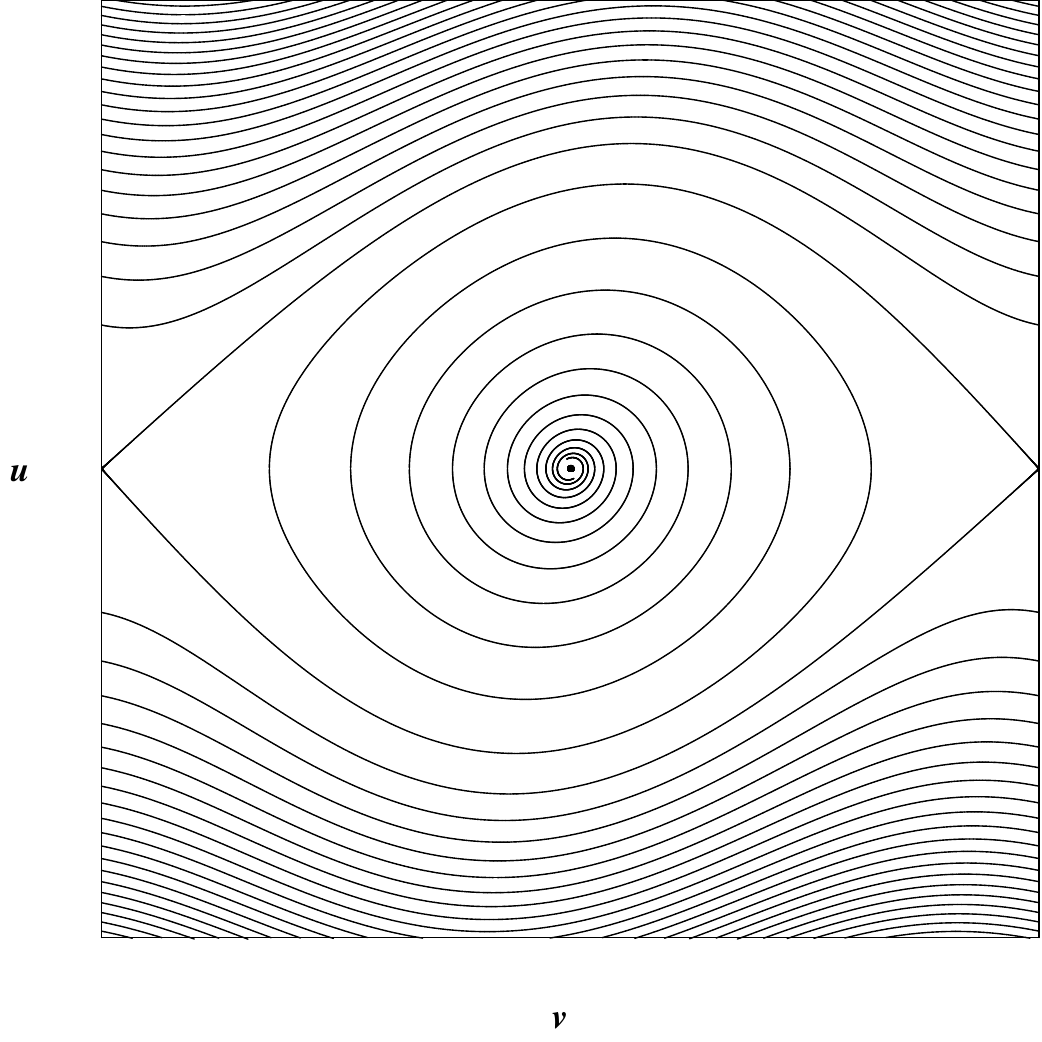}&
			\includegraphics[width=140pt,height=140pt]{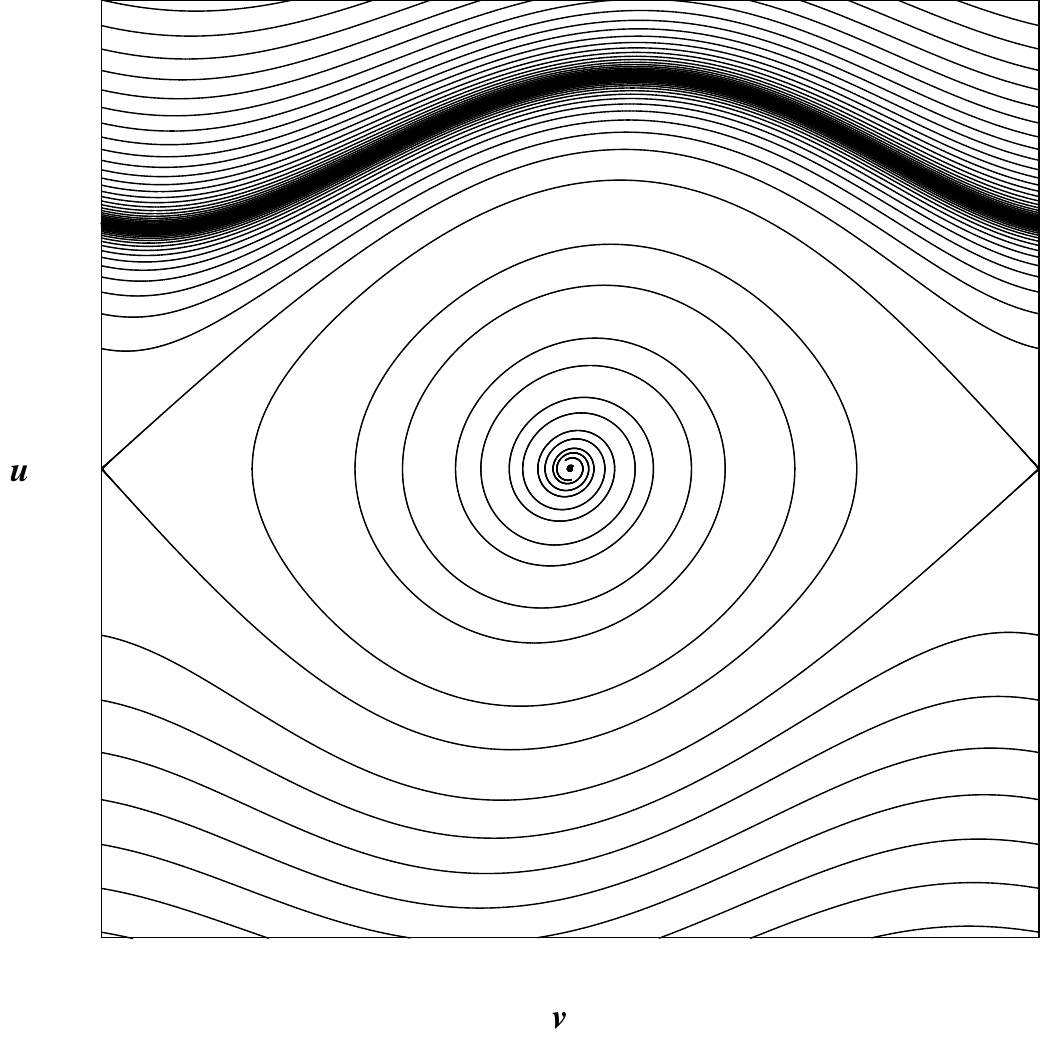}\\
			\footnotesize{(a) $B_2(I_{nm_1m_2})=0$}& \footnotesize{(b) $B_2(I_{nm_1m_2})\neq 0$}\\
			\includegraphics[width=140pt,height=140pt]{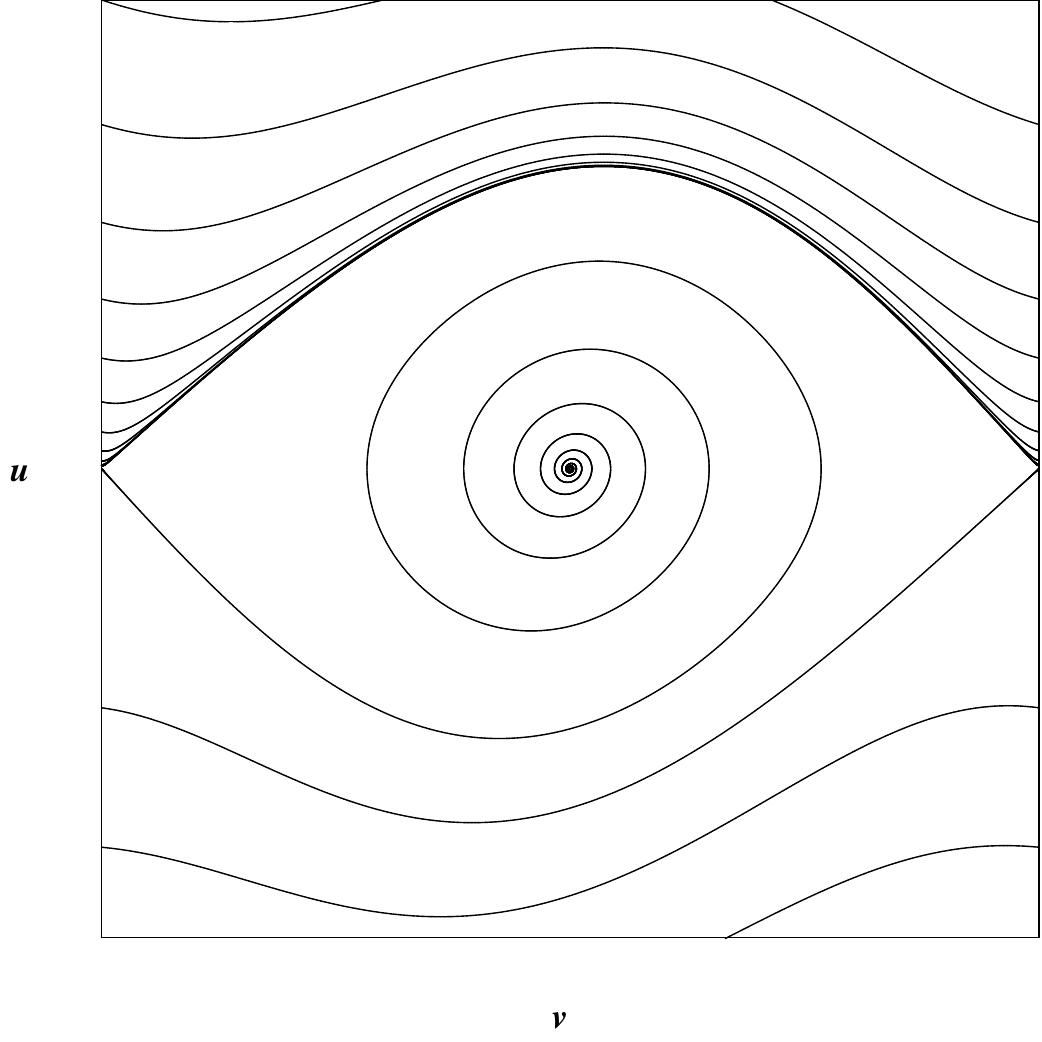}&
			\includegraphics[width=140pt,height=140pt]{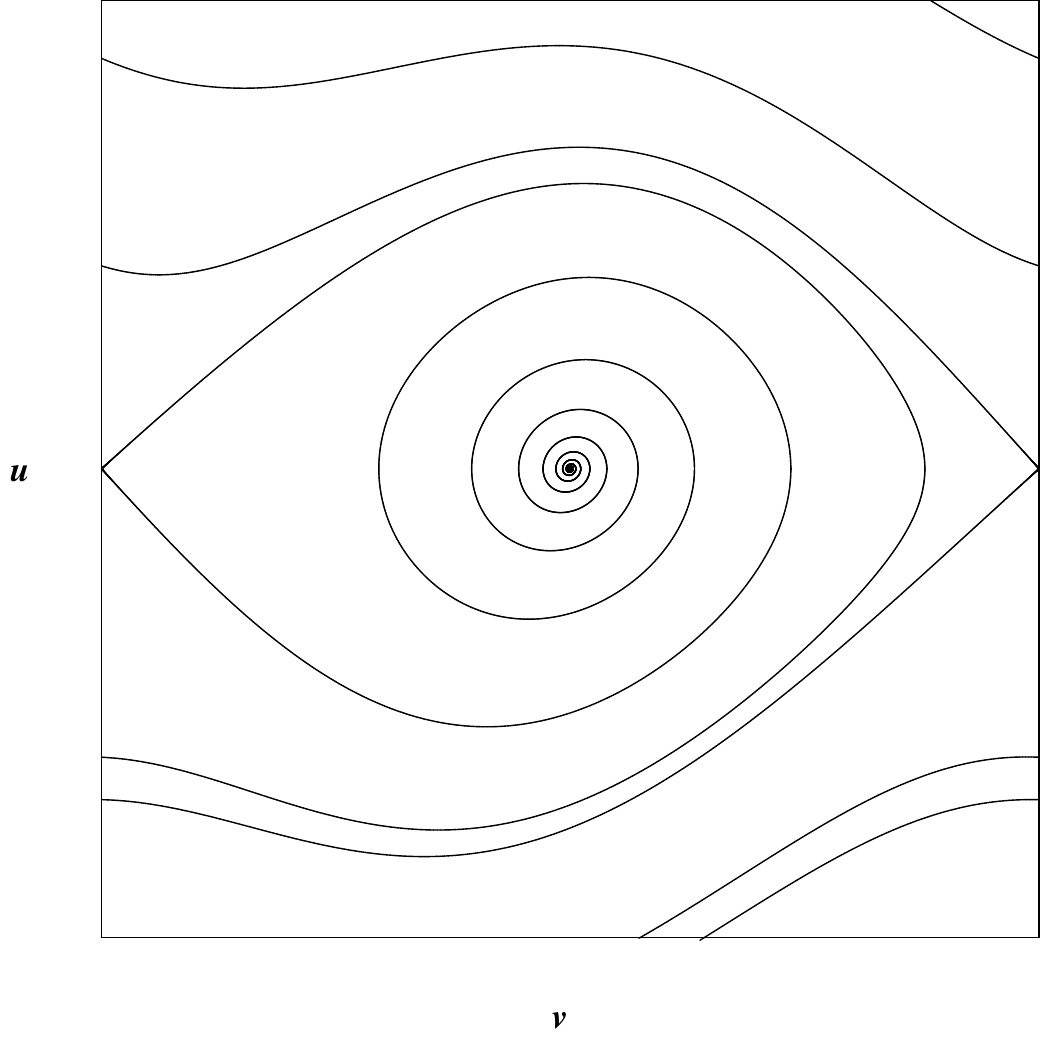}\\
			\footnotesize{(c) $B_2(I_{nm_1m_2})\neq 0$}& \footnotesize{(d) $B_2(I_{nm_1m_2})\neq 0$}\\
		\end{tabular}
	\end{center}
	\renewcommand{\figurename}{Fig}
	\caption{Possible phase portraits of system~\eqref{eq13} for an alternating function $\widetilde{\sigma}(v;I_{nm_1m_2})$.}\label{fig2}
\end{figure}

Let us consider the case of an alternating function $\widetilde{\sigma}(v;I_{nm_1m_2})$.  Possible phase portraits of system~\eqref{eq13} on the phase cylinder $\{v(mod(2\pi/n)),u\}$ in this case are shown in Fig~\ref{fig2}.
The structure of impassable resonance zone is shown in Figs~\ref{fig2}(a) and (b): the phase point tends to a stable equilibrium state at $\tau\to\infty$ for any initial point (a); for any initial point on the lower half-cylinder and on an unstable limit cycle located on the upper half-cylinder (b). The bifurcation case, when the limit cycle merges into the separatrix loop on the upper half-cylinder, is shown in Fig~\ref{fig2}(c). The question of limit cycles can be solved using the Poincar\'e--Pontryagin generating function, similarly to how it is done for a perturbed autonomous system~\eqref{eq5}. The structure of partly passable resonance zone is shown in Fig~\ref{fig2}(d). 

\begin{figure}[h!]
	\begin{center}
		\begin{tabular}{ccc}
			\includegraphics[width=140pt,height=140pt]{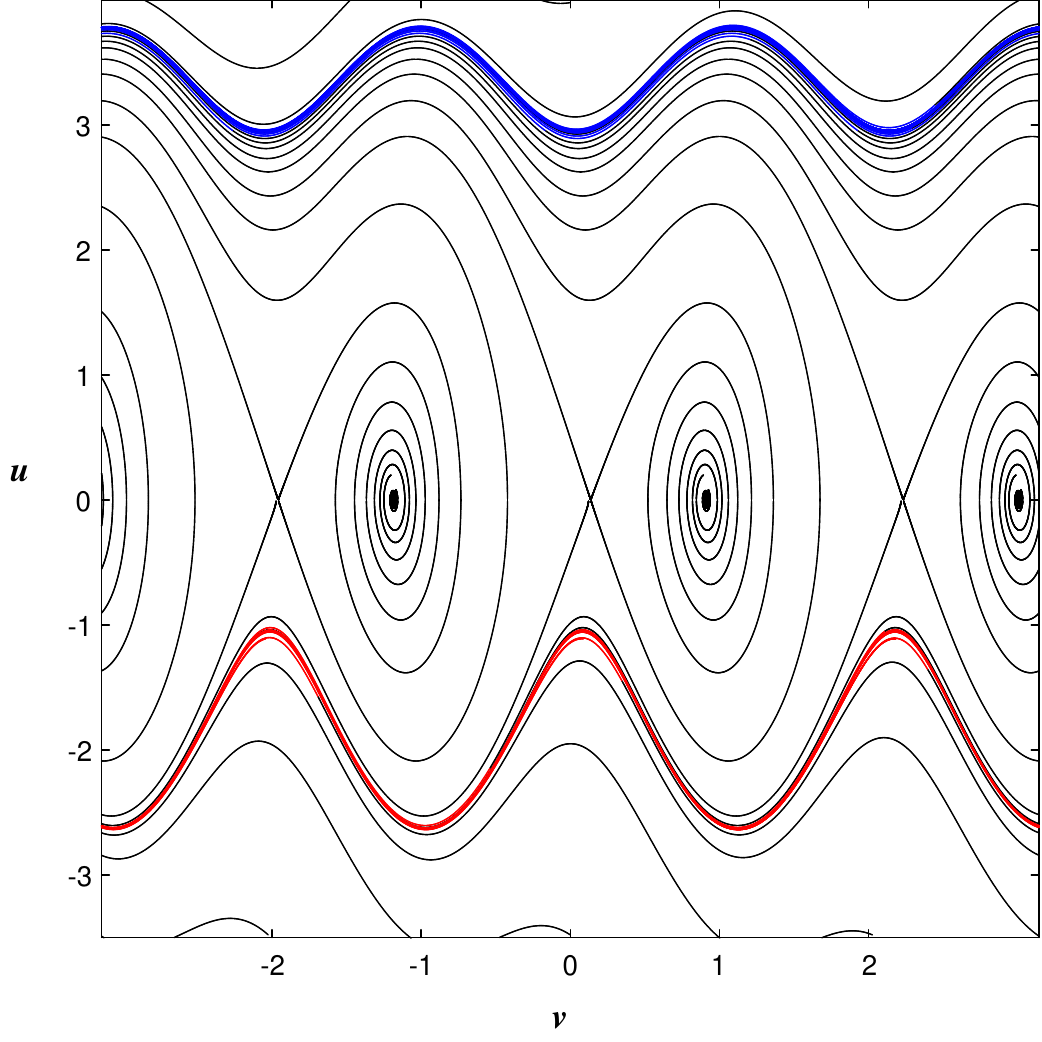}&
			\includegraphics[width=140pt,height=140pt]{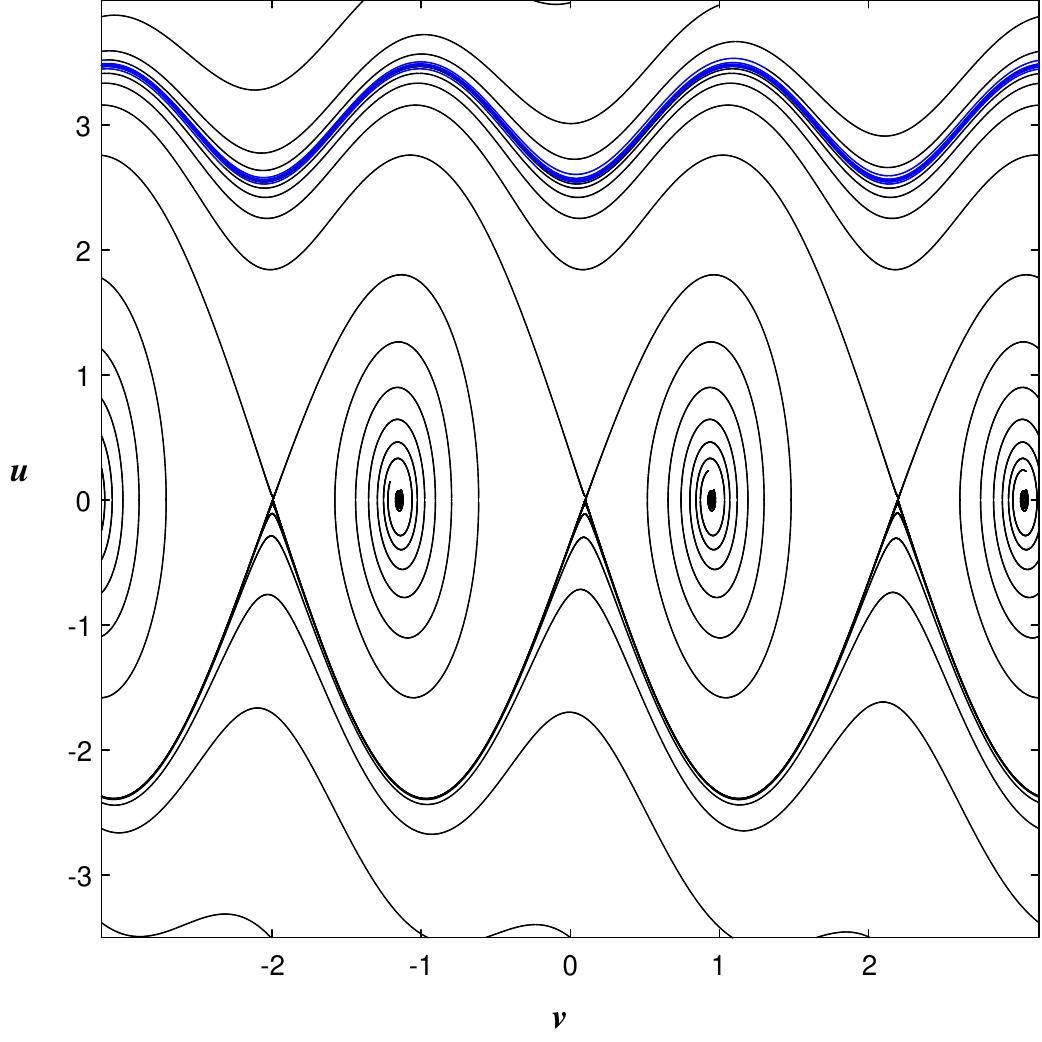}&
			\includegraphics[width=140pt,height=140pt]{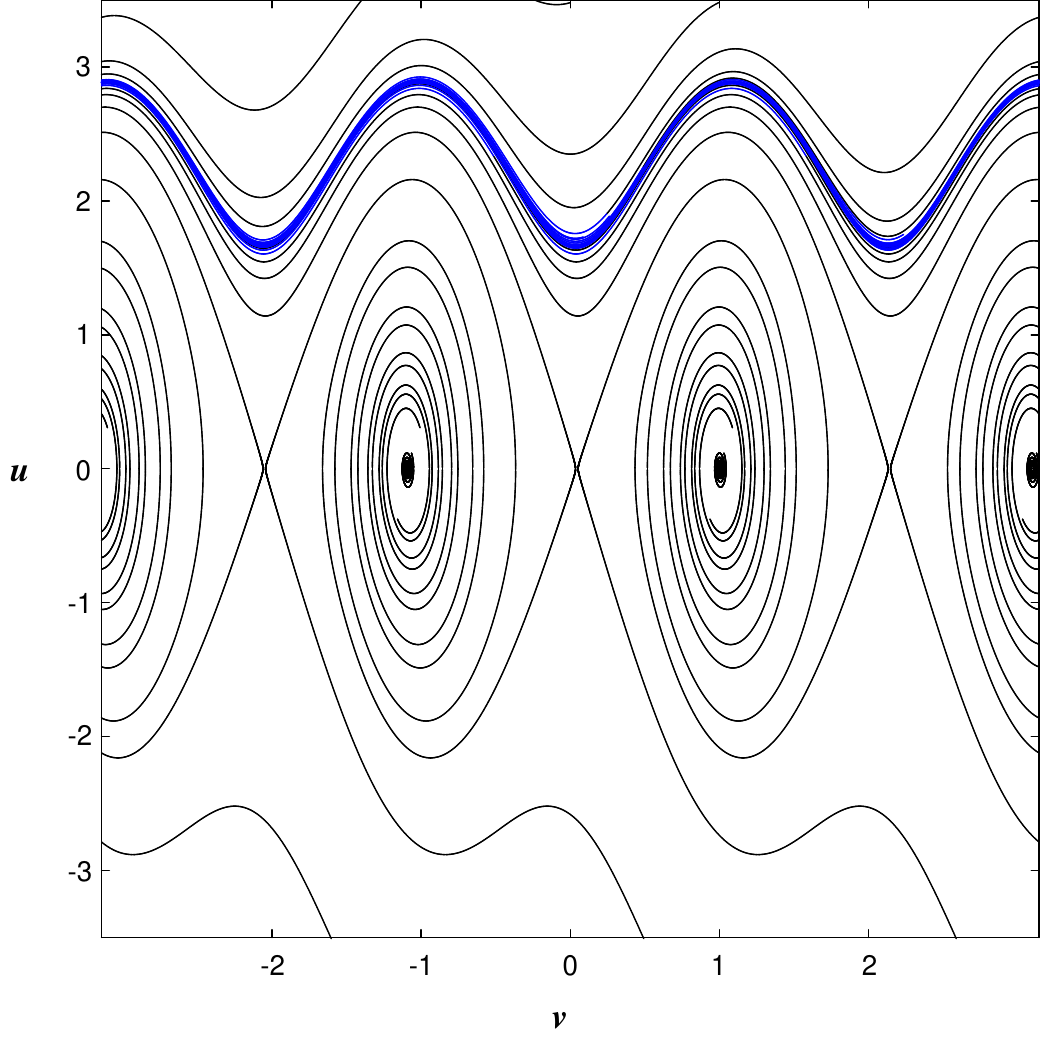}\\
			\footnotesize{(a)}& \footnotesize{(b)}& \footnotesize{(c)}\\
			\includegraphics[width=140pt,height=140pt]{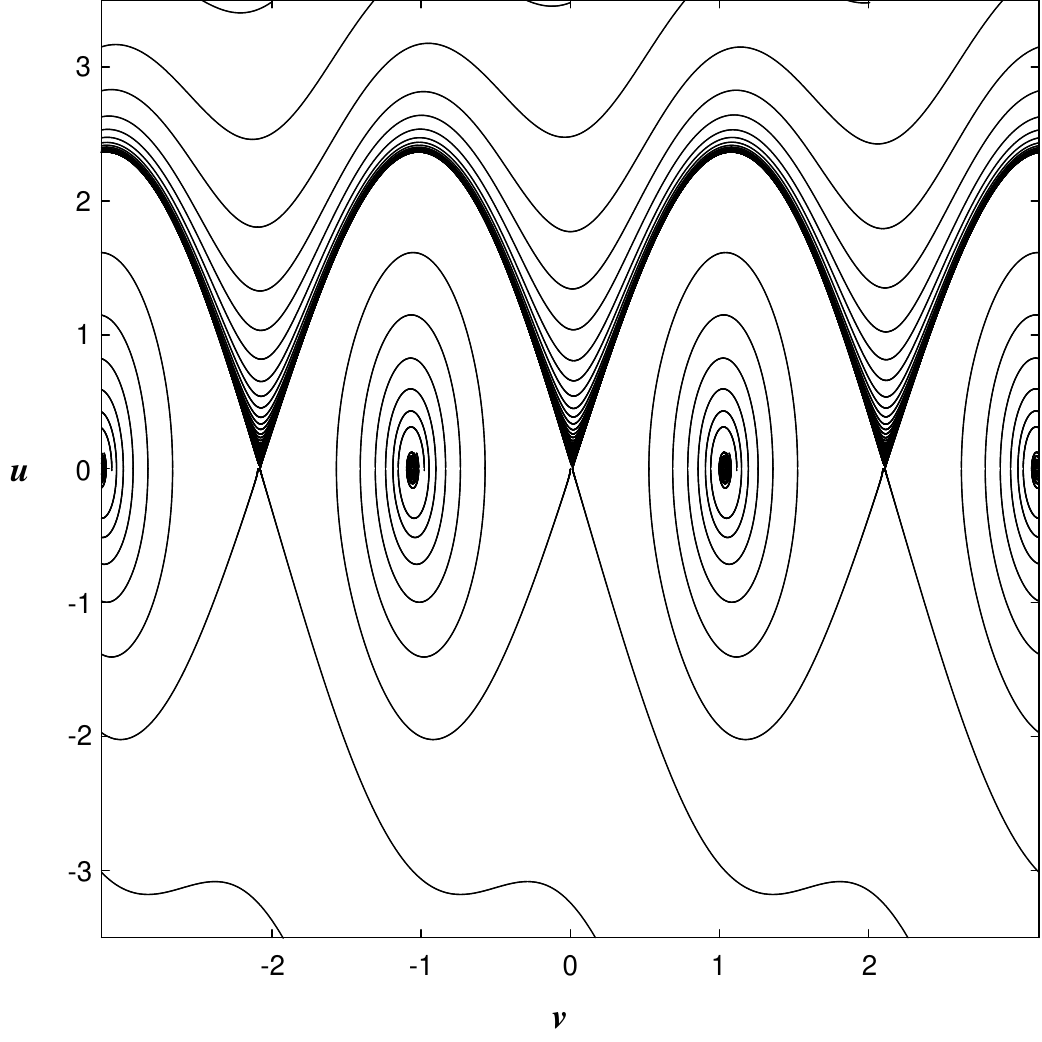}&
			\includegraphics[width=140pt,height=140pt]{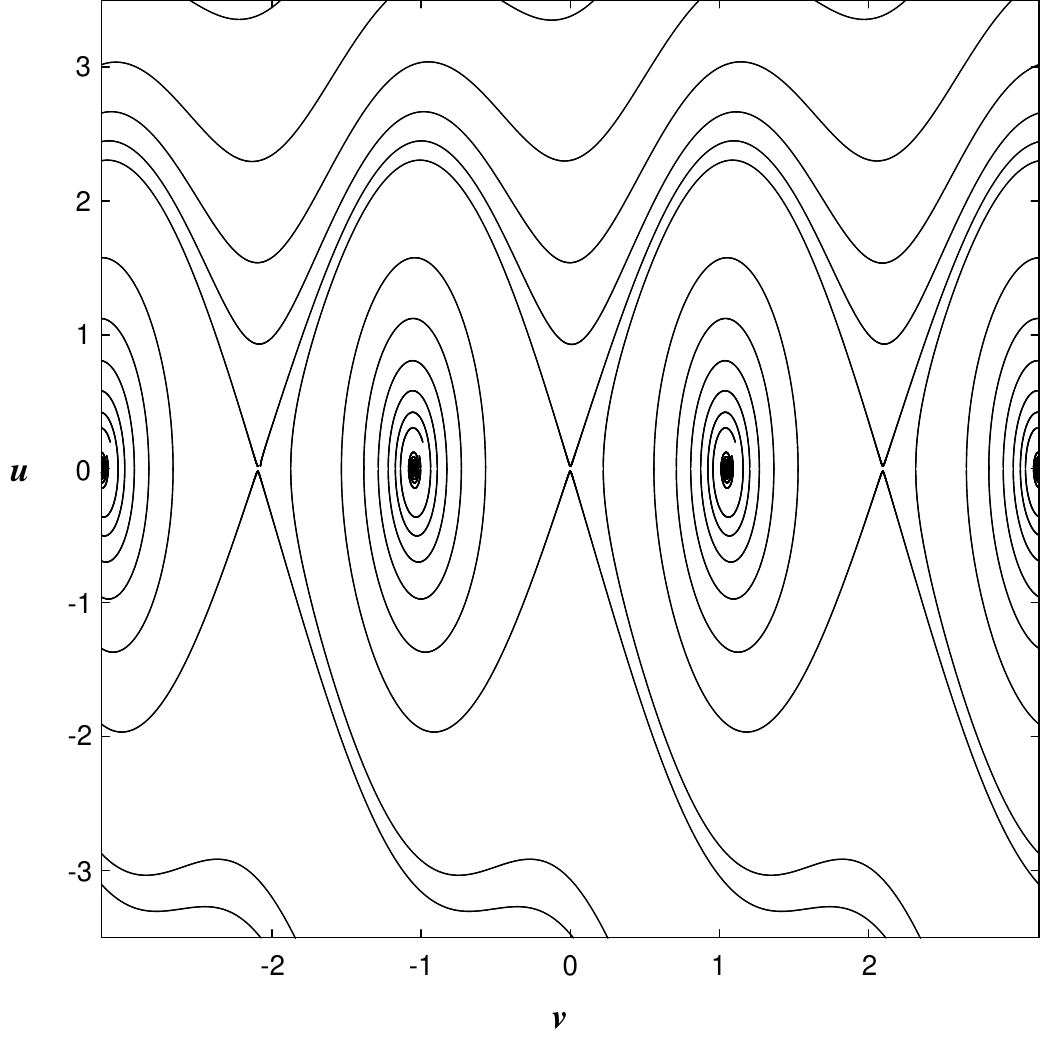}\\
			\footnotesize{(d)}& \footnotesize{(e)}\\
		\end{tabular}
	\end{center}
	\renewcommand{\figurename}{Fig}
	\caption{Possible phase portraits of system~\eqref{eq10} for $n=3$ in the case of a constant sign function $\sigma\neq 0$.}\label{fig3}
\end{figure}

We assign $B_0(I_{nm_1m_2})=\mu \gamma_1$, $B_1(I_{nm_1m_2}) = \mu \gamma_2$. For nonzero parameters $\gamma_1$ and $\gamma_2$, possible phase portraits of system~\eqref{eq10} for $n=3$ on the phase cylinder $\{v(mod(2\pi)),u\}$ in the case of a constant sign function $\sigma\neq 0$ are shown in Fig~\ref{fig3}. The focus will be on the phase cylinder of system~\eqref{eq10} instead of the center, and there will also be an asymmetry in the formation of limit cycles on the upper and lower half-cylinders of different stability (shown in different colors). Therefore, a bifurcation case is possible when one of the limit cycles (on the upper or lower phase half cylinder) merges into the corresponding separatrix loop; as the parameter changes further, it disappears.

\begin{figure}[h!]
	\begin{center}
		\begin{tabular}{ccc}
			\includegraphics[width=140pt,height=140pt]{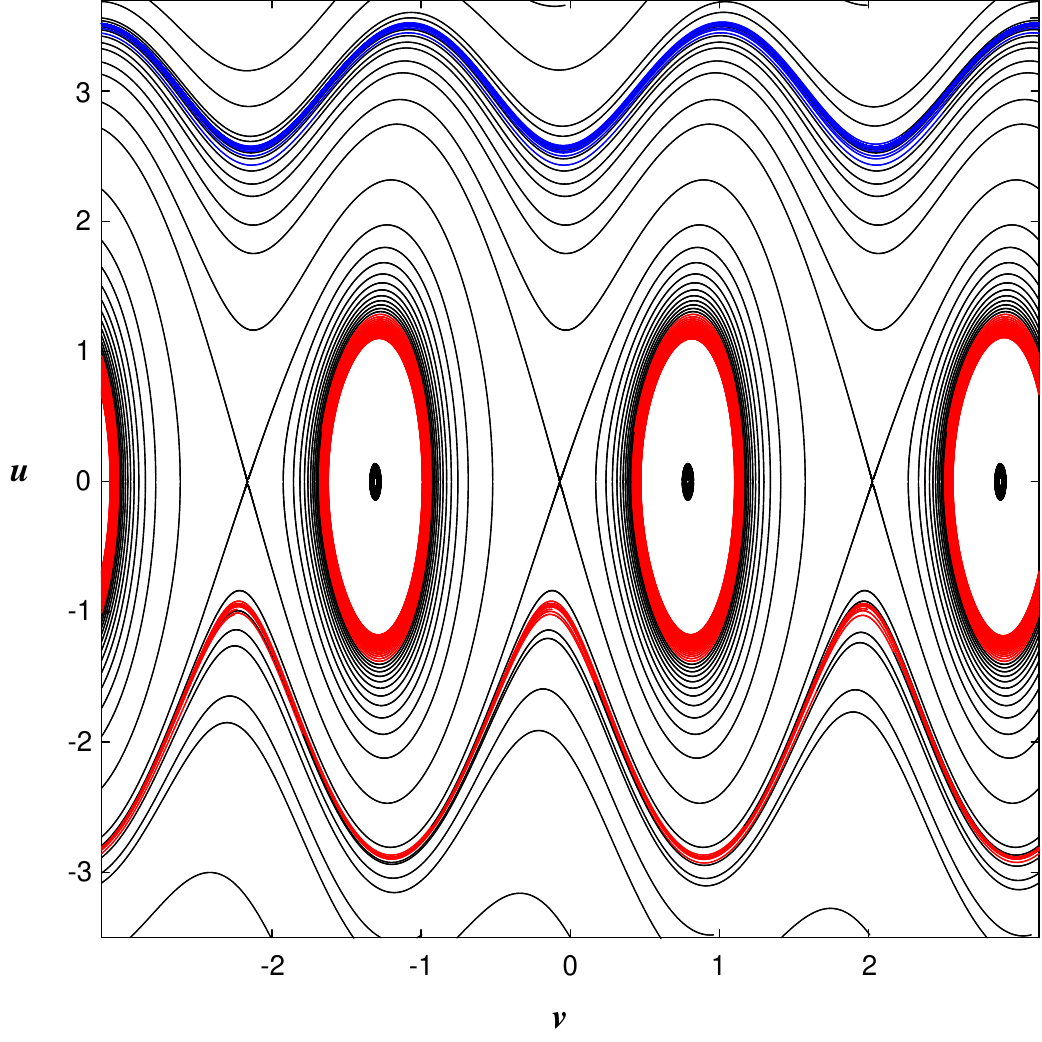}&
			\includegraphics[width=140pt,height=140pt]{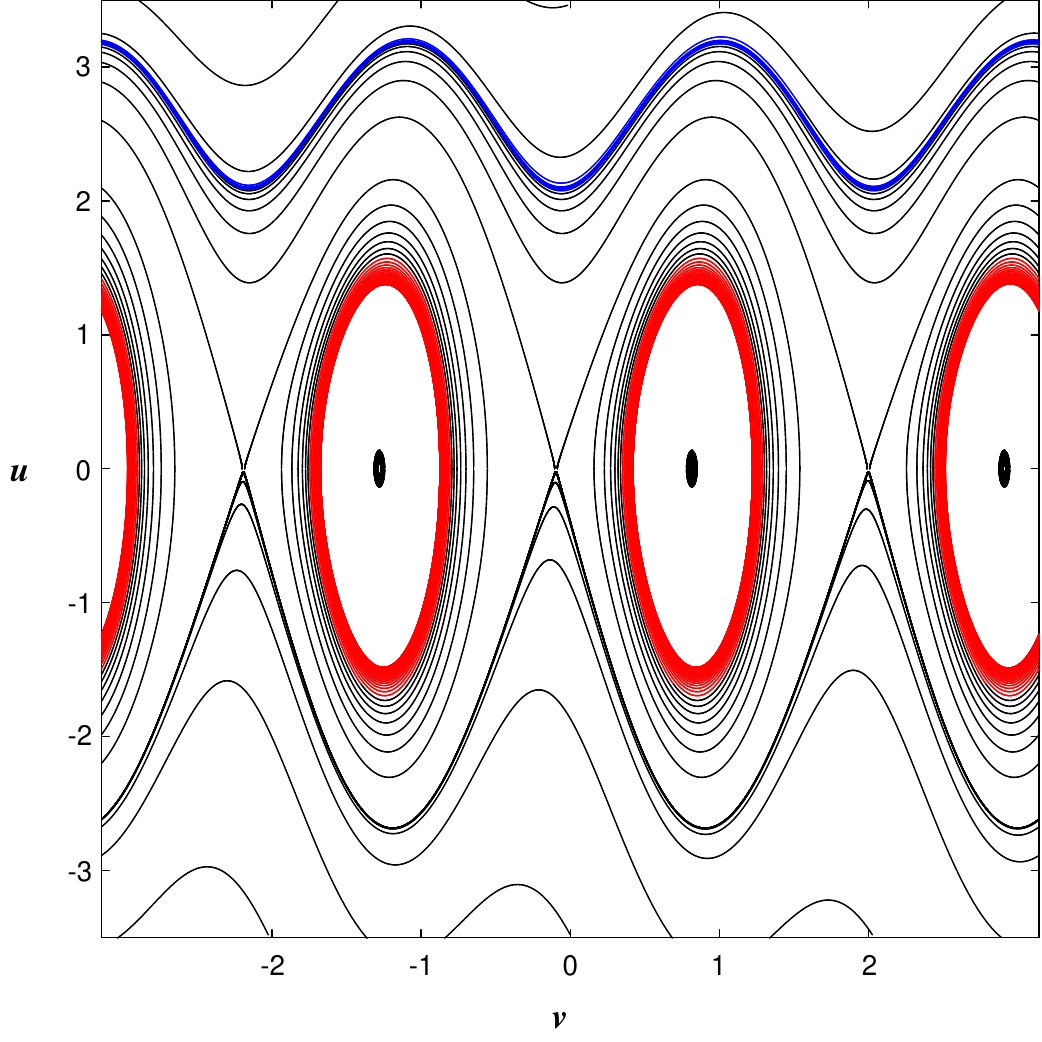}&
			\includegraphics[width=140pt,height=140pt]{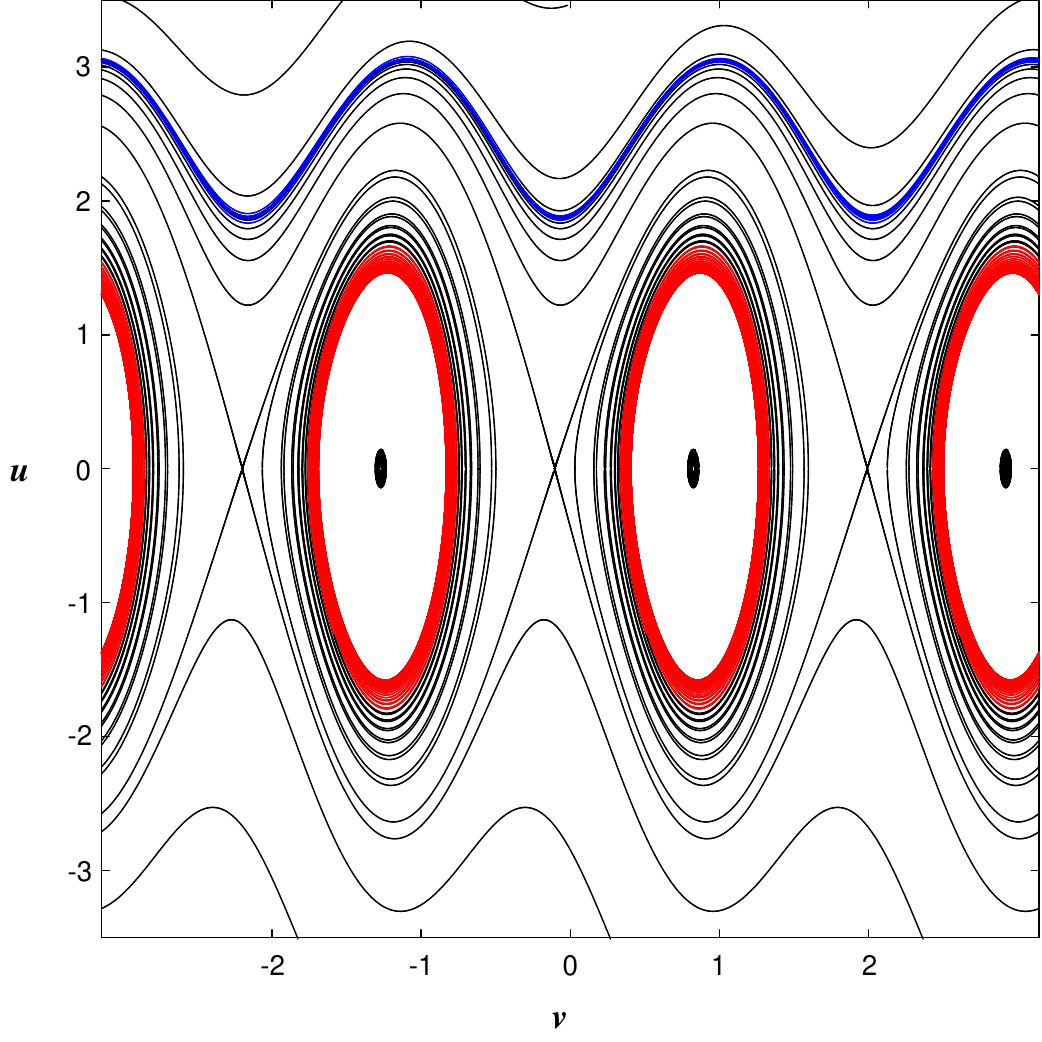}\\
			\footnotesize{(a)}& \footnotesize{(b)}& \footnotesize{(c)}\\
			\includegraphics[width=140pt,height=140pt]{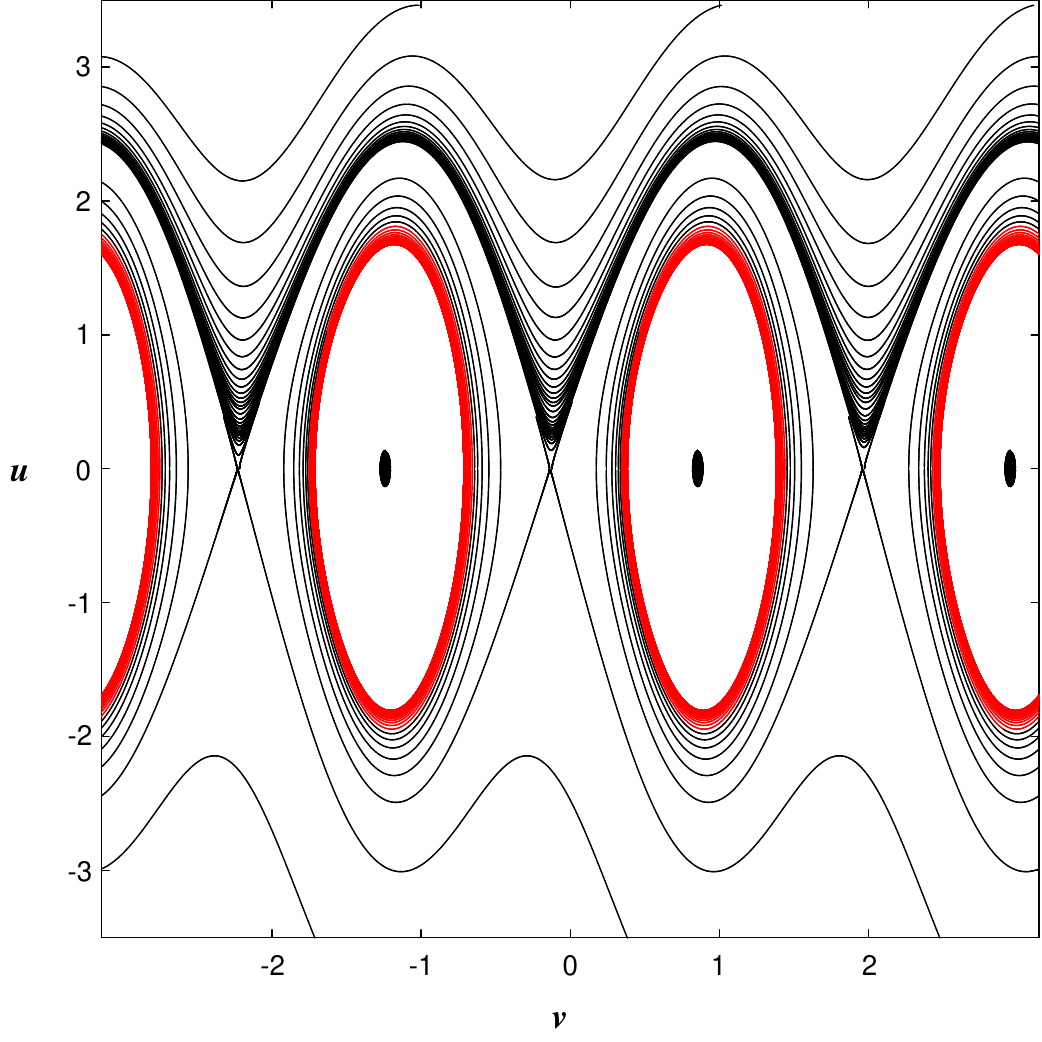}&
			\includegraphics[width=140pt,height=140pt]{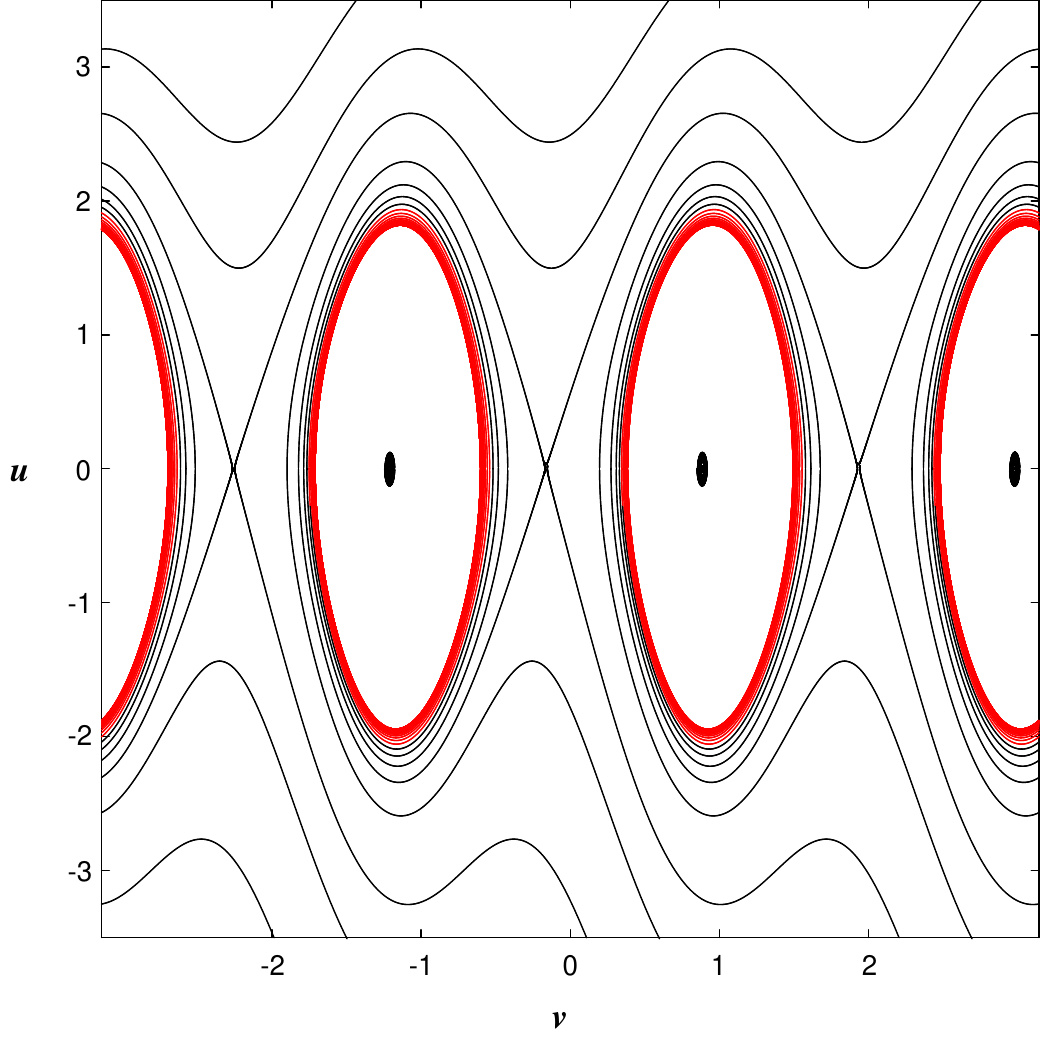}\\
			\footnotesize{(d)}& \footnotesize{(e)}\\
		\end{tabular}
	\end{center}
	\renewcommand{\figurename}{Fig}
	\caption{Possible phase portraits of system~\eqref{eq10} for $n=3$ in the case of an alternating function $\sigma$.}\label{fig4}
\end{figure}

In the case of an alternating function $\sigma$, the existence of a limit cycle in the region of oscillatory motions is possible. Possible phase portraits of system~\eqref{eq10} for $n=3$ on the phase cylinder $\{v(mod(2\pi)),u\}$ in the case of an alternating function $\sigma$ are shown in Fig~\ref{fig4}.

Cases of synchronization of oscillations are shown in Figs~\ref{fig3}(c) and ~\ref{fig4}(c). Figures~\ref{fig3}(a), \ref{fig3}(e), \ref{fig4}(a) and \ref{fig4}(e) show the structures of partly passable resonance zones. Bifurcation cases are shown in Figs~\ref{fig3}(b), \ref{fig3}(d), \ref{fig4}(b) and \ref{fig4}(d). 

In~\cite{30}, the phenomenon of synchronization of oscillations was called restricted. This is due to the possibility of the existence of a limit cycle in the region of rotational motion of the pendulum; therefore, not all phase points in the resonance zone tend to a stable periodic regime, which corresponds to a stable equilibrium state of the averaged system.

\section{Example for illustration}

Let us consider the equation~\eqref{eqExample}. The autonomous equation ($p_2=0$ and $p_3=0$) was studied in the paper~\cite{31}. It was found that the Poincar\'e--Pontryagin generating function for region of oscillatory motion (filled with closed phase curves $\dot{x}^2/2-\cos{x}=h$, $h\in (-1,1)$) of the unperturbed pendulum has the form
\begin{equation}\label{eq17}
\begin{split}
B_0=B_0(k(h))&=\frac{8}{105\pi}\{(1-k^2)(105+(128k^4-80k^2+3)p_1){\bf K}(k)+\\
&+(-105+(2k^2-1)(128k^4-128k^2+3)p_1){\bf E}(k)\},
\end{split}
\end{equation}
where ${\bf K}(k)$ and ${\bf E}(k)$ are complete elliptic integrals of the first and second kind respectively; $k=\sqrt{(1+h)/2}\in (0,1)$ is their module. In ~\cite{31}, we found the bifurcation value $p_1\approx -8.481$, for which there is a double limit cycle. In this case, the Poincar\'e--Pontryagin generating equation $B_0(k)=0$ has a double root  $k=k_*\approx 0.759$.

Now we choose such $\omega_1$ and $\omega_2$ so that the level $k=k_*$ coincides with the resonance level $k=k_{311}$. From the resonance condition $3\omega(k_*)=\omega_1+\omega_2$, where $\omega(k_*)=\displaystyle\frac{\pi}{2{\bf K}(k_*)}$, we find, for example, $\omega_1=1$ and $\omega_2\approx1.448$. 

Let us calculate system ~\eqref{eq10} for our case. This needs to be calculated:
\begin{equation}\label{eq18}
\widetilde{A}(v)=\frac 1{36\pi^2}\int\limits_0^{6\pi}\int\limits_0^{6\pi}(p_2xy+p_3)\frac{\partial x}{\partial \theta}\cos\theta_1\sin\theta_2\,d\theta_1\,d\theta_2 ,
\end{equation}
\begin{equation}\label{eq19}
\widetilde{P_0}(v)=\frac 1{36\pi^2}\int\limits_0^{6\pi}\int\limits_0^{6\pi}\frac{\partial}{\partial I}\left((p_2xy+p_3)\frac{\partial x}{\partial \theta}\right)\cos\theta_1\sin\theta_2 \,d\theta_1\,d\theta_2 ,
\end{equation}
\begin{equation}\label{eq20}
Q_0(v)=-\frac 1{36\pi^2}\int\limits_0^{6\pi}\int\limits_0^{6\pi}[(-1+p_1\cos{3x}+p_2x\cos\theta_1\sin\theta_2)y+p_3\cos\theta_1\sin\theta_2] \frac{\partial x}{\partial I}\,d\theta_1\,d\theta_2 ,
\end{equation}
\begin{equation}\label{eq21}
\widetilde{P_1}(v)=\frac 1{72\pi^2}\int\limits_0^{6\pi}\int\limits_0^{6\pi}\frac{\partial^2}{\partial I^2}\left((p_2xy+p_3)\frac{\partial x}{\partial \theta}\right)\cos\theta_1\sin\theta_2 \,d\theta_1\,d\theta_2 ,
\end{equation}
\begin{equation}\label{eq22}
Q_1(v)=-\frac 1{36\pi^2}\int\limits_0^{6\pi}\int\limits_0^{6\pi}\frac{\partial}{\partial I}\left([(-1+p_1\cos{3x}+p_2x\cos\theta_1\sin\theta_2)y+p_3\cos\theta_1\sin\theta_2]\frac{\partial x}{\partial I}\right)\,d\theta_1\,d\theta_2 ,
\end{equation}
\begin{equation}\label{eq23}
B_1(k)=\frac{dB_0(k)}{dI}=\frac{dB_0(k)}{dk}\frac{dk}{dh}\frac{dh}{dI}=\frac{\pi}{8k{\bf K}(k)}\frac{dB_0(k)}{dk}, \ B_2(k)=\frac{d^2B_0(k)}{2dI^2}=\frac{\pi}{16k{\bf K}(k)}\frac{dB_1(k)}{dk},
\end{equation}
\begin{equation}\label{eq24}
\widetilde{\sigma}(v)=\frac 1{36\pi^2}\int\limits_0^{6\pi}\int\limits_0^{6\pi}p_2x\cos\theta_1\sin\theta_2 \,d\theta_1\,d\theta_2 ,
\end{equation}
\begin{equation}\label{eq25}
b_1=\frac{\pi}{8k{\bf K}(k)}\frac{d\omega(k)}{dk}, \  b_2=\frac{\pi}{16k{\bf K}(k)}\frac{db_1(k)}{dk}, \ b_3=\frac{\pi}{24k{\bf K}(k)}\frac{db_2(k)}{dk},
\end{equation}
where $x(k,\theta)=2\arcsin{(k\sn(2{\bf K}\theta/\pi,k))}$, $y(k,\theta)=2k\cn(2{\bf K}\theta/\pi,k)$, $\theta=v+\displaystyle\frac{\theta_1+\theta_2}{3}$ ($\sn(u,k)$ is the Jacobi elliptic sine, $\cn(u,k)$ is the Jacobi elliptic cosine). Here and below we put $k=k_*=k_{311}$.

In calculating integrals~\eqref{eq18}, \eqref{eq19} and \eqref{eq21}, using the well-known Fourier series expansion for the Jacobi elliptic cosine~\cite{33}, we obtain $$\displaystyle\frac{\partial x}{\partial \theta}=\displaystyle\frac{y}{\omega}=\displaystyle\frac{4k{\bf K}}{\pi}\cn(2{\bf K}\theta/\pi,k)=8\sum_{j=1}^{\infty}\displaystyle\frac{a^{j-1/2}}{1+a^{2j-1}}\cos(2j-1)\theta ,$$
$$\displaystyle\frac{\partial }{\partial I}\left(\displaystyle\frac{\partial x}{\partial \theta}\right)=\displaystyle\frac{\pi^3}{4k^2(1-k^2){\bf K}^3}\sum_{j=1}^{\infty}\displaystyle\frac{(2j-1)a^{j-1/2}\left(1-a^{2j-1}\right)}{\left(1+a^{2j-1}\right)^2}\cos(2j-1)\theta ,$$
\begin{equation*}
\begin{split}
\frac{\partial^2 }{\partial I^2}\left(\frac{\partial x}{\partial \theta}\right)&=\frac{\pi^4}{32k^4(1-k^2)^2{\bf K}^5}\left(((1+k^2){\bf K}-3{\bf E})\sum_{j=1}^{\infty}\frac{(2j-1)a^{j-1/2}\left(1-a^{2j-1}\right)}{\left(1+a^{2j-1}\right)^2}\cos(2j-1)\theta\right.+\\
&+\left.\frac{\pi^2}{4{\bf K}}\sum_{j=1}^{\infty}\frac{(2j-1)^2a^{j-1/2}\left(1-6a^{2j-1}+a^{2(2j-1)}\right)}{\left(1+a^{2j-1}\right)^3}\cos(2j-1)\theta\right) ,
\end{split}
\end{equation*}
where $a=exp\left(-\pi\displaystyle\frac{{\bf K}(\sqrt{1-k^2})}{{\bf K}(k)} \right)$. In calculating integrals~\eqref{eq20} and \eqref{eq22}, expansions for the functions $\displaystyle\frac{\partial x}{\partial I}=\displaystyle\frac{\pi\left(\sn(.)\dn(.)-\cn(.)zn(.)\right)}{4k(1-k^2){\bf K}}$ (here the dot denotes the arguments $(2{\bf K}\theta/\pi,k)$; $\dn(u,k)$ is the delta amplitude, $\zn(u,k)$ is the Jacobi zeta function) and $\displaystyle\frac{\partial }{\partial I}\left(\displaystyle\frac{\partial x}{\partial I}\right)$ in a Fourier series
$$\displaystyle\frac{\partial x}{\partial I}=\displaystyle\frac{\pi^3}{4k^2(1-k^2){\bf K}^3}\sum_{j=1}^{\infty}\displaystyle\frac{a^{j-1/2}\left(1-a^{2j-1}\right)}{\left(1+a^{2j-1}\right)^2}\sin(2j-1)\theta ,$$
\begin{equation*}
\begin{split}
\frac{\partial }{\partial I}\left(\frac{\partial x}{\partial I}\right)&=\frac{\pi^4}{32k^4(1-k^2)^2{\bf K}^5}\left(((1+k^2){\bf K}-3{\bf E})\sum_{j=1}^{\infty}\frac{a^{j-1/2}\left(1-a^{2j-1}\right)}{\left(1+a^{2j-1}\right)^2}\sin(2j-1)\theta\right.+\\
&+\left.\frac{\pi^2}{4{\bf K}}\sum_{j=1}^{\infty}\frac{(2j-1)a^{j-1/2}\left(1-6a^{2j-1}+a^{2(2j-1)}\right)}{\left(1+a^{2j-1}\right)^3}\sin(2j-1)\theta\right)
\end{split}
\end{equation*}
were obtained. When calculating the integral \eqref{eq24}, we have $$x=8\sum_{j=1}^{\infty}\displaystyle\frac{a^{j-1/2}}{(2j-1)\left(1+a^{2j-1}\right)}\sin(2j-1)\theta .$$

Thus, calculating \eqref{eq18} -- \eqref{eq25}, we obtain an averaged system of the following form
\begin{equation}\label{eq26}
\left\{\begin{aligned}
\dot{u} &=p_3\widetilde{A_1}\sin{3v}+p_2\widetilde{A_2}\cos{3v}+B_0+\mu\left(p_2\widetilde{\sigma}\cos{3v}+B_1\right)u+\\&+\mu^2\left[\left(p_3\left(\widetilde{P_{11}}-\frac{3b_2}{b_1}Q_{01}\right)\sin{3v}\right.\right.+\left.\left.p_2\left(\widetilde{P_{12}}+\frac{3b_2}{b_1}Q_{02}\right)\cos{3v}+B_2\right)u^2\right.-\\&-\left.\frac{1}{b_1}\left(p_3\widetilde{P_{01}}\sin{3v}+p_2\widetilde{P_{02}}\cos{3v}\right)\left(p_3Q_{01}\cos{3v}+p_2Q_{02}\sin{3v}\right)\right],\\
\dot{v} &=b_1u+\mu b_2u^2+\mu^2\left[b_3u^3+\left(p_3\left(Q_{11}-\frac{2b_2}{b_1}Q_{01}\right)\cos{3v}+p_2\left(Q_{12}-\frac{2b_2}{b_1}Q_{02}\right)\sin{3v}\right)u\ \right],
\end{aligned}\right.
\end{equation}
where $$\widetilde{A_1}=-2\displaystyle\frac{a^{3/2}}{1+a^3} , \ \widetilde{A_2}=\displaystyle\frac{1}{4\pi\omega}\int_{0}^{2\pi}xy^2\sin{3\theta}d\theta ,$$
$$\widetilde{P_{01}}=\displaystyle\frac{\pi^3}{16k^2(1-k^2){\bf K}^3}\displaystyle\frac{3a^{3/2}(a^3-1)}{(1+a^3)^2} , \ \widetilde{P_{02}}=\displaystyle\frac{1}{4\pi\omega}\int_{0}^{2\pi}\displaystyle\frac{\partial }{\partial I}\left(xy^2\right)\sin{3\theta}d\theta ,$$
$$\widetilde{P_{11}}=\displaystyle\frac{\pi^4}{1024k^4(1-k^2)^2{\bf K}^6}\displaystyle\frac{3a^{3/2}}{(1+a^3)^3}\left(4((1+k^2){\bf K}-3{\bf E}){\bf K}(a^6-1)-3\pi^2(a^6-6a^3+1)\right) ,$$
$$\widetilde{P_{12}}=\displaystyle\frac{1}{8\pi\omega}\int_{0}^{2\pi}\displaystyle\frac{\partial^2 }{\partial I^2}\left(xy^2\right)\sin{3\theta}d\theta , \ Q_{01}=\displaystyle\frac{\widetilde{P_{01}}}{3} , \ Q_{02}=\displaystyle\frac{1}{4\pi}\int_{0}^{2\pi}xy\displaystyle\frac{\partial x}{\partial I}\cos{3\theta}d\theta ,$$
$$Q_{11}=\displaystyle\frac{2\widetilde{P_{11}}}{3} , \ Q_{12}=\displaystyle\frac{1}{4\pi}\int_{0}^{2\pi}\displaystyle\frac{\partial}{\partial I}\left(xy\right)\displaystyle\frac{\partial x}{\partial I}\cos{3\theta}d\theta , \ \widetilde{\sigma}=\displaystyle\frac{2a^{3/2}}{3(1+a^3)} ,$$
$$B_1=-\displaystyle\frac{1}{15}\left[15+\left(128k^4-144k^2+31\right)p_1-\left(256k^4-256k^2+46\right)p_1\displaystyle\frac{{\bf E}}{{\bf K}}\right] ,$$
\begin{equation*}
\begin{split}
B_2&=\displaystyle\frac{p_1\pi}{120k^2(1-k^2){\bf K}^3}\left[\left(384k^6-656k^4+295k^2-23\right){\bf K}^2-\left(768k^6-1280k^4+558k^2-46\right){\bf K}{\bf E}\right.-\\
&-\left.\left(128k^4-128k^2+23\right){\bf E}^2\right] ,
\end{split}
\end{equation*}
$$b_1=\frac{\pi^2}{16}\frac{\left(k^2-1\right){\bf K}+{\bf E}}{k^2\left(k^2-1\right){\bf K}^3} , \ b_2=-\frac{\pi^3}{256}\frac{\left(k^2-1\right){\bf K}^2-2\left(k^2-2\right){\bf K}{\bf E}-3{\bf E}^2}{k^4\left(k^2-1\right)^2{\bf K}^5} ,$$
$$b_3=\frac{\pi^4}{6144}\frac{\left(k^4+2k^2-3\right){\bf K}^3-\left(2k^4+3k^2-13\right){\bf K}^2{\bf E}+5\left(k^2-5\right){\bf K}{\bf E}^2+15{\bf E}^3}{k^6\left(k^2-1\right)^3{\bf K}^7} ,$$
$B_0$ is determined by formula \eqref{eq17}.

Possible phase portraits of the averaged system \eqref{eq26} obtained using the WInSet software~\cite{13},~\cite{14} are shown in Figures ~\ref{fig1}--~\ref{fig3}.


\section{Conclusion}

The problem of the action of time-quasiperiodic perturbations on a self-oscillating system plays an important role in the theory of oscillations. As in the case of time-periodic perturbations, it leads to the phenomenon of synchronization of oscillations, first studied by Andronov and Witt~\cite{36} on the example of the Van der Pol equation in the quasi-linear case. Later, this phenomenon was studied by Morozov and Shilnikov~\cite{37} for systems close to arbitrary two-dimensional Hamiltonian ones (with time-periodic perturbation).

The existence of quasi-periodic solutions in systems close to Hamiltonian ones has been considered in many papers. However, studies of the effect of quasi-periodic perturbations on the self-oscillating system appeared for the first time, apparently, in the papers of Morozov and his students~\cite{1,3,7,8,11}.

In this paper, the problem of the effect of quasi-periodic perturbations on systems close to arbitrary two-dimensional Hamiltonian in the case when perturbed autonomous systems have a double limit cycle was studied. Averaged systems are obtained for the case when the resonance level is close to the level from which a double limit cycle appears under the influence of an autonomous perturbation. These systems are not amenable to investigation in the general case.
This means the need to consider new examples and study the behavior of solutions of these averaged systems in each case.

\section{Acknowledgments}
This research was supported by the Russian Foundation for Basic Research (project No. 18-01-00306), the Russian Science Foundation (project No. 19-11-00280), and the Ministry of Science and Higher Education of Russian Federation under grant No. 0729-2020-0036.

	\renewcommand{\refname}{{\rm


\centerline{REFERENCES}}}
	
	\begin{thebibliography}{99}
		
	\bibitem{1} Morozov A.D., Dragunov T.N. On quasi-periodic perturbations of Duffing equation, {\it Discontinuity, Nonlinearity, and Complexity}, 2016, vol.\,5, no. 4, pp. 377--386. https://doi.org/10.5890/DNC.2016.12.005	
	
	\bibitem{2} Kuznetsov A.P., Kuznetsov S.P., Sedova Y.V. Pendulum system with an infinite number of equilibrium states and quasiperiodic dynamics, {\it Nelineinaya Dinamika}, 2016, vol. 12, no. 2, pp. 223--234. (In Russian). https://doi.org/10.20537/nd1602005	
	
	\bibitem{3}  Morozov A.D., Morozov K.E. Quasiperiodic perturbations of two-dimensional Hamiltonian systems, {\it Differential Equations}, 2017, vol. 53, no. 12, pp. 1557--1566. https://doi.org/10.1134/S0012266117120047
	
	\bibitem{4} Stankevich N.V., Kuznetsov A.P., Popova E.S., Seleznev E.P. Experimental diagnostics of multi-frequency quasiperiodic oscillations, {\it Communications in Nonlinear Science and Numerical Simulation}, 2017, vol. 43, pp. 200--210. https://doi.org/10.1016/j.cnsns.2016.07.007
	
	\bibitem{5} Truong T.Q., Tsubone T., Sekikawa M., Inaba N. Complicated quasiperiodic oscillations and chaos from driven piecewise-constant circuit: chenciner bubbles do not necessarily occur via simple phase-locking, {\it Physica D}, 2017, vol. 341, pp. 1--9. https://doi.org/10.1016/j.physd.2016.09.008
	
	\bibitem{6} Jiang T., Yang, Z., Jing, Z. Bifurcations and Chaos in the Duffing Equation with Parametric Excitation and Single External Forcing, {\it International Journal of Bifurcation and Chaos}, 2017, vol. 27, no. 08, p. 1750125-1-31. https://doi.org/10.1142/S0218127417501255
	
	\bibitem{7} Morozov A.D., Morozov K.E. On synchronization of quasiperiodic oscillations, {\it Russian Journal of Nonlinear Dynamics}, 2018, vol. 14, no. 3, pp. 367--376. https://doi.org/10.20537/nd180307
	
	\bibitem{8} Morozov A.D., Morozov K.E. Global dynamics of systems close to Hamiltonian ones under nonconservative quasi-periodic perturbation, {\it Russian Journal of Nonlinear Dynamics}, 2019, vol. 15, no. 2, pp. 187--198. https://doi.org/10.20537/nd190208
	
	\bibitem{9} Kuznetsov A.P., Kuznetsov S.P., Shchegoleva N.A., Stankevich N.V. Dynamics of coupled generators of quasiperiodic oscillations: Different types of synchronization and other phenomena, {\it Physica D}, 2019, vol. 398, pp. 1--12. https://doi.org/10.1016/j.physd.2019.05.014
	
	\bibitem{10} Guan Y., Gupta V., Wan M., Li L.K.B. Forced synchronization of quasiperiodic oscillations in a thermoacoustic system, {\it Journal of Fluid Mechanics}, 2019, vol. 879, pp. 390--421. https://doi.org/10.1017/jfm.2019.680
	
	\bibitem{11} Morozov A.D., Morozov K.E. On quasi-periodic parametric perturbations of Hamiltonian systems, {\it Russian Journal of Nonlinear Dynamics}, 2020, vol. 16, no. 2, pp. 369--378. https://doi.org/10.20537/nd200210
	
	\bibitem{12} Kolos M., Shahzadi M., Stuchlik Z. Quasi-periodic oscillations around Kerr-MOG black holes, {\it The European Physical Journal C}, 2020, vol. 80, no. 133, pp. 1--12. https://doi.org/10.1140/epjc/s10052-020-7692-5	
	
	\bibitem{13} Sato M., Hyodo H., Biwa T., Delage R. Synchronization of thermoacoustic quasiperiodic oscillation by periodic external force, {\it Chaos}, 2020, vol. 30, no. 6, p. 063130. https://doi.org/10.1063/5.0004381
	
	\bibitem{14} Blekhman, I.I. {\it Synchronization in science and technology}, New York: ASME Press, 1988, 255 p. 	
	
	\bibitem{15} Pikovsky A., Rosenblum M., Kurths J. {\it Synchronization: A universal concept in nonlinear science}, England, Cambridge: Cambridge University Press, 2001. https://doi.org/10.1119/1.1475332
	
	\bibitem{16} Kuznetsov A.P., Sataev I.R., Stankevich N.V., Turukina L.V. {\it Physics of quasiperiodic oscillations}, Saratov: Publishing center ``Science'', 2013, 252 p. (In Russian).
	
	\bibitem{17} Anishchenko V., Nikolaev S., Kurths J. Winding number locking on a two-dimensional torus: Synchronization of quasiperiodic motions, {\it Phys. Rev. E.}, 2006, vol. 73, p. 056202. https://doi.org/10.1103/PhysRevE.73.056202
	
	\bibitem{18} Kuznetsov A.P., Kuznetsov S.P., Stankevich N.V. A simple autonomous quasiperiodic self-oscillator, {\it Communications in Nonlinear Science and Numerical Simulation}, 2010, vol. 15, pp. 1676--1681. https://doi.org/10.1016/j.cnsns.2009.06.027
	
	\bibitem{19} Kuznetsov A.P., Kuanetsov S.P., Mosekilde E., Stankevich N.V. Generators of quasiperiodic oscillations with three-dimensional phase space, {\it European Physical Journal Special Topics}, 2013, vol. 222, pp. 2391--2398. https://doi.org/10.1140/epjst/e2013-02023-x
	
	\bibitem{20} Berger M.S., Chen  Y. Y. Forced quasiperiodic and almost periodic oscillations of nonlinear Duffing equations, {\it Nonlinear Analysis: Theory, Methods and Applications}, 1992, vol. 19, no. 3, pp. 249--257. https://doi.org/10.1016/0362-546X(92)90143-3
	
	\bibitem{21} Liu B., You J. Quasiperiodic solutions of Duffing’s Equations, {\it Nonlinear Analysis: Theory, Methods \& Applications}, 1998, vol. 33, no. 6, p. 645--655. https://doi.org/10.1016/S0362-546X(98)00662-2	
	
	\bibitem{22} Grischenko A.D., Vavriv D.M., Dynamics of pendulum with a quasiperiodic perturbation, {\it Technical Physics}, 1997, vol. 42, no. 10, pp. 1115--1120. https://doi.org/10.1134/1.1258787	
	
	\bibitem{23} Wang R.Q., Deng J., Jing Z.J. Chaos control in duffing system, {\it Chaos, Solitons and Fractals}, 2006, vol. 27, no. 1, pp. 249--257. https://doi.org/10.1016/j.chaos.2005.03.038
	
	\bibitem{24} Jing Z., Yang Z., Jiang T. Complex dynamics in Duffing-van der Pol equation, {\it Chaos, Solitons and Fractals}, 2006, vol. 27, no. 3, pp. 722--747. https://doi.org/10.1016/j.chaos.2005.04.044
	
	\bibitem{25} Jing Z.J., Huang J.C., Deng J. Complex dynamics in three-well Duffing system with two external forcings, {\it Chaos, Solitons and Fractals}, 2007, vol. 33, no. 3, pp.795--812. https://doi.org/10.1016/j.chaos.2006.03.071
	
	\bibitem{26} Ravichandran V., Chinnathambi V, Rajasekar S. Homoclinic bifurcation and chaos in Duffing oscillator driven by an amplitude-modulated force, {\it Physica A: Statistical Mechanics and its Applications}, 2007, vol. 376, pp. 223--236. https://doi.org/10.1016/j.physa.2006.11.003		
	
	\bibitem{27} Hale J.K. {\it Oscillations in nonlinear systems}, Mineola, New York: Dover Publications, Inc., 2015, 192 p.
	
	\bibitem{28} Morozov A.D. {\it Quasi-conservative systems: cycles, resonances and chaos}, Singapore: World Sci, in ser. Nonlinear Science, ser. A, V. 30, 1998, 340 p. https://doi.org/10.1142/3238
	
	\bibitem{29} Morozov A.D. {\it Resonance, cycles and chaos in quasi-conservative systems}, Izhevsk: Regular and Chaotic Dynamics, Institute of Computer Science, 2005, 424 p. (In Russian).
	
	\bibitem{30} Morozov A.D., Mamedov E.A. On a double cycle and resonances, {\it Regular and Chaotic Dynamics}, 2012, vol. 17, no. 1, pp. 63--71. https://doi.org/10.1134/S1560354712010066
	
	\bibitem{31} Kostromina O.S. On limit cycles, resonance and homoclinic structures in asymmetric pendulum-type equation, {\it Vestn. Udmurtsk. Univ. Mat. Mekh. Komp. Nauki}, 2019, vol. 29, no. 2, pp. 228--244. https://doi.org/10.20537/vm190207
	
	\bibitem{32} Mel'nikov V.K. On stability of a center under periodic in time perturbations, {\it Works of Moscow Math. Socity}, 1963, vol. 12, pp. 3--53. (In Russian). http://mi.mathnet.ru/eng/mmo137
	
	\bibitem{33} Gradshteyn I.S., Ryzhik I.M. Table of integrals, series, and products, {\it Academic Press}, 2007, 1171 p.
	
	\bibitem{34} Morozov A.D., Dragunov T.N., Boykova S.A., Malysheva O.V. {\it Invariant sets for Windows}, World Scientific, 1999, 272 p. https://doi.org/10.1142/4220
	
	\bibitem{35} Morozov A.D., Dragunov T.N. {\it Visualization and analysis of invariant sets of dynamical systems}, Moscow--Izhevsk: Publishing House of the Institute of Computer Research, 2003, 303 p. (In Russian).
	
	\bibitem{36} Andronov A.A., Witt A.A. Zur theorie des mitnehmens von van der Pol, {\it Archiv fur Elektrotechnik}, 1930, vol. 24, no. 1, pp. 99--110. https://doi.org/10.1007/BF01659580
	
	\bibitem{37} Morozov A.D., Shil'nikov L.P. On nonconservative periodic systems close to two-dimensional Hamiltonian, {\it Journal of Applied Mathematics and Mechanics}, 1983, vol. 47, issue 3, pp. 327--334. https://doi.org/10.1016/0021-8928(83)90058-8
		
\end{thebibliography}
\end{document}